\documentclass[12pt,reqno]{amsart}
\pdfoutput=1
\usepackage{amsmath,amssymb,amscd,graphicx,soul}
\usepackage[all]{xy}
\usepackage[latin1]{inputenc}
\usepackage[T1]{fontenc}
\usepackage{hyperref}
\usepackage{ifthen}
\usepackage{fullpage}
\usepackage{subsupscripts}
\usepackage{xcolor}
\usepackage{faktor}
\usepackage{tikz-cd}

\newboolean{workmode}
\setboolean{workmode}{false}

\newboolean{sections}
\setboolean{sections}{true}

\input xy
\xyoption{all}

\makeatletter
\def\th@plain{%
  \thm@notefont{}
  \itshape 
}
\def\th@definition{%
  \thm@notefont{}
  \normalfont 
}
\makeatother

\newcommand{\calG}{{\mathcal{G}}}
\newcommand{\calH}{{\mathcal{H}}}
\newcommand{\calI}{{\mathcal{I}}}

\newcommand{\calK}{{\mathcal{K}}}
\newcommand{\calL}{{\mathcal{L}}}
\newcommand{\calM}{{\mathcal{M}}}
\newcommand{\calN}{{\mathcal{N}}}

\newcommand{\calP}{{\mathcal{P}}}
\newcommand{\calQ}{{\mathcal{Q}}}

\newcommand{\id}{{\operatorname{id}}}  
\newcommand{\ID}{{\operatorname{ID}}} 
\newcommand{\im}{{\operatorname{im}}}  
\newcommand{\inv}{{\operatorname{inv}}} 
\newcommand{\LieGpoid}{{\mathbf{LieGpoid}}} 
\newcommand{\pr}{{\operatorname{pr}}} 
\newcommand{\PR}{{\operatorname{PR}}} 
\newcommand{\ic}{{\operatorname{inc}}} 
\newcommand{\Stab}{{\operatorname{Stab}}} 

\newcommand{\G}{{\operatorname{G}}} 

\newcommand{\eps}{\varepsilon}  
\newcommand{\ftimes}[2]{{\lrsubscripts{\times}{#1}{#2}}} 
\newcommand{\toto}{{~\rightrightarrows~}} 
\newcommand{\ot}{{~\overset{\simeq}{\longleftarrow}~}} 
\newcommand{\wtimes}[2]{{\lrsubscripts{\overset{\raisebox{-1pt}{\tiny$\mathrm{w}$}}{\times}}{#1}{#2}}} 

\makeatletter
\DeclareRobustCommand*{\mfaktor}[3][]
{
   { \mathpalette{\mfaktor@impl@}{{#1}{#2}{#3}} }
}
\newcommand*{\mfaktor@impl@}[2]{\mfaktor@impl#1#2}
\newcommand*{\mfaktor@impl}[4]{
   \settoheight{\faktor@zaehlerhoehe}{\ensuremath{#1#2{#3}}}%
   \settoheight{\faktor@nennerhoehe}{\ensuremath{#1#2{#4}}}%
      \raisebox{-0.5\faktor@zaehlerhoehe}{\ensuremath{#1#2{#3}}}%
      \mkern-4mu\diagdown\mkern-5mu%
      \raisebox{0.5\faktor@nennerhoehe}{\ensuremath{#1#2{#4}}}%
}
\makeatother
 
\newcommand{\lfaktor}[2]{{\mfaktor{#2}{#1}}}

\newcommand{\ifwork}[1]{\ifthenelse{\boolean{workmode}}{#1}{}}
\newcommand{\comment}[1]{}
\newcommand{\mute}[1]{}
\newcommand{\printname}[1]{}

\ifwork{
\renewcommand{\comment}[1]{{\marginpar{*}\ \scriptsize{#1}\ }}
\renewcommand{\printname}[1]
    {{\color{brown}{\makebox[0pt]{\hspace{-1.0in}\raisebox{8pt}{\tiny #1}}}}}
}

\newcommand{\ifsection}[2]{\ifthenelse{\boolean{sections}}{#1}{#2}}

\theoremstyle{plain}
\ifsection{
    \newtheorem{theorem}{Theorem}[section]
}
{
    \newtheorem{theorem}{Theorem}
}

\newtheorem{proposition}[theorem]{Proposition}

\newtheorem{corollary}[theorem]{Corollary}
\newtheorem{lemma}[theorem]{Lemma}

\theoremstyle{definition}
\newtheorem{definition}[theorem]{Definition}
\newtheorem{example}[theorem]{Example}
\newtheorem{remark}[theorem]{Remark}

\definecolor{jaw}{rgb}{0,.5,0}

\definecolor{carla}{rgb}{0, 0, 1} 
\definecolor{laura}{rgb}{.4, 0, .6} 

\title{Bicategories of Action Groupoids}
\date{\today}

\author{Carla Farsi}
\address{Department of Mathematics, University of Colorado at Boulder, Boulder, Colorado 80309}
\email{carla.farsi@colorado.edu}

\author{Laura Scull}
\address{Department of Mathematics, Fort Lewis College, Durango, Colorado 81301}
\email{scull\_l@fortlewis.edu}

\author{Jordan Watts}
\address{Department of Mathematics, Central Michigan University, Pearce Hall 214, Mount Pleasant, MI 48859, USA}
\email{jordan.watts@cmich.edu}

\newcommand{\eqactgpd}{{\mathbf{ActGpd}}} 
\newcommand{\ana}{{\mathbf{AnaLieGpoid}}}
\newcommand{\gana}{{\mathbf{AnaActGpd}}}
\newcommand{\ssw}{{\sf{sW}}}
\newcommand{\we}{{\sf{W}}}

\newcommand{\ES}[1]{{\textup {\bf ES}_{#1}}} 

\newcommand{\Ff}[1]{{\textup {\bf FF}_{#1}}} 
\DeclareRobustCommand\longtwoheadrightarrow{\relbar\joinrel\twoheadrightarrow} 

\begin{document}

\begin{abstract}
We prove that the 2-category of action Lie groupoids localised in the following three different ways yield equivalent bicategories: localising at equivariant weak equivalences \`a la Pronk, localising using surjective submersive equivariant weak equivalences and anafunctors \`a la Roberts, and localising at all weak equivalences.  These constructions generalise the known case of representable orbifold groupoids.  We also show that any weak equivalence between action Lie groupoids is isomorphic to the composition of two particularly nice forms of equivariant weak equivalences.
\end{abstract}



\maketitle

\section{Introduction}\label{s:intro}

The study of  topological or Lie groupoids up to Morita equivalence, generated by a topologised/differential geometric version of  categorical weak equivalence, appears in many contexts.    One standard way to do so, as in Pronk \cite{pronk},  is to define a bicategory of fractions of Lie groupoids $\LieGpoid[\we^{-1}]$  in which all of the Morita equivalences become invertible $1$-cells. In this bicategory, a $1$-cell from $\calG $ to $ \calH$ is defined by a span of functors $\calG \leftarrow \calK \to \calH$ where the functor $\calK \to \calG$ is a weak equivalence, and a $2$-cell is defined as an equivalence class of diagrams with natural transformations.   There is also a ``smaller'' localisation using anafunctors \`a la Roberts \cite{roberts2012, roberts2021}. By additionally requiring the  weak equivalences used in the span to be surjective submersive on objects, we obtain a bicategory with similar $1$-cells, but whose $2$-cells are actual natural transformations chosen to represent  the equivalence classes defining the $2$-cells of Pronk.     Roberts proves that this bicategory is equivalent to the bicategory of fractions of Pronk. 
Understanding conditions under which a class of $1$-cells of a general bicategory admits a localisation is a current field of study;  see also \cite{AV,pronk-scull:bicategories}.  The localisation of topological groupoids is made explicit in \cite{C-D-W}, where the authors discuss Lie groupoids but details are only provided for the topological case.  In \cite{watts:bicateg}, a detailed development of localisations of diffeological groupoids using Pronk's and Roberts' approaches is presented, comparing and connecting them to the theory of bibundles in the diffeological groupoid context established in \cite{vdS}, as well as to stacks over diffeological spaces.

An important class of groupoids comes from Lie group actions.   Having an action groupoid allows for the application of an equivariant  functor such as cohomology.   The question then becomes whether the invariant  respects Morita equivalence.  Functoriality ensures that our chosen invariant is unchanged under \emph{equivariant} weak equivalences, but \emph{a priori} there is no mechanism for checking more general weak equivalences.  Thus we wish to know whether Morita equivalent action groupoids are also Morita equivalent via equivariant weak equivalences.
More precisely:   if we consider the full sub-$2$-category of action groupoids,  can we  create a bicategory of fractions inverting the equivariant weak equivalences?  If so, is this equivalent to the full sub-bicategory we get by inverting more  general weak equivalences and then restricting to action groupoids? 

If the answer to both questions is yes,  we have a mechanism  for transferring equivariant techniques to groupoids: choose a groupoid that is Morita equivalent to an action groupoid, apply some equivariant functor, and check that the result is Morita invariant.  
This strategy has been successfully applied,  
for example, to orbifold groupoids  Morita equivalent to an action groupoid:   in \cite[Proposition 5.13]{pronk-scull} they define Bredon cohomology under a mild condition on the coefficient systems.

In this paper, instead of focusing on orbifolds, we generalise our group actions to those which have any subset of the following properties: free, locally free, transitive, effective, compact, discrete, proper, and being Morita equivalent to a proper \'etale Lie groupoid.  (There is, of course, much redundancy when multiple properties are taken together.)  Given action groupoids satisfying any selection of these properties, we localise at equivariant weak equivalences, and we show that the resulting bicategory is equivalent to the localisation at all weak equivalences; see Theorem \ref{t:eqvt Lie gpds}.   This yields affirmative answers to the two questions above for these action groupoids. Moreover, we show that we can construct this localisation using the method of Roberts, giving us a smaller and more concrete category to work with; see Proposition~\ref{p:eqvt bicat 2}.   This is a first step towards generalising \cite{pronk-scull} to more general action groupoids, and defining Bredon cohomology as a Morita invariant in a more general context.
As an additional step towards understanding when equivariant functors might be Morita invariant, we also show that the  decomposition of equivariant weak equivalences used in \cite[Proposition 3.5]{pronk-scull} also applies in our  more general setting.   This allows us to break down equivariant weak equivalences into two specific types: projections and inclusions.   This decomposition has proved useful in other contexts such as topological complexity \cite{ACGO}.  

The paper is structured as follows:  
Section~\ref{s:bicateg frac} contains background information on Lie groupoids, and Section~\ref{s:localisation} background on the localisation of $\LieGpoid$ at the class of weak equivalences, as well as at surjective submersive weak equivalences.   Sections~\ref{s:eqactgpd} and \ref{s:eqvt-to-full} contains our main results about action groupoids, answering the two questions posed above.  In particular, in Section~\ref{s:eqactgpd}, we construct the localisation of action groupoids satisfying any subset of properties from a fixed list (see \eqref{e:P} and Proposition~\ref{p:eqvt bicat}) and show that  equivariant weak
equivalences with these chosen properties can be decomposed into an equivariant projection and an equivariant inclusion  (see Theorem \ref{t:decomposition}).    In Section~\ref{s:eqvt-to-full} we show that this is equivalent to the full sub-bicategory of such action groupoids (see Theorem~\ref{t:eqvt Lie gpds}).  

This paper is written with a wide audience in mind, and so some details (for instance, proofs of smoothness of maps) are spelled out. Throughout, we refer to \cite{benabou,JY,nlab} for categorical definitions such as bicategory, pseudofunctor, $2$-category, and $2$-commutative, and all of their constituent parts.

 \textbf{Acknowledgements:}  We thank Dorette Pronk for pointing us to certain references on bicategories and $2$-categories, and the anonymous referee for many helpful suggestions. Carla Farsi was partially supported by the Simons Foundation Collaboration Grant \#523991.

\section{The 2-Category of Lie Groupoids}\label{s:bicateg frac}
In this section we discuss the 2-category of Lie groupoids and its properties, with special attention to the notion of weak equivalence, which gives rise to the ubiquitous notion of Morita equivalence.  

We begin by clarifying our setting.  Throughout this paper, ``smooth'' means infinitely-differentiable, and all manifolds are smooth and without boundary.  By a ``curve'', we mean a smooth map $p\colon I:=(-\eps,\eps)\to M$ where $M$ is a smooth manifold and $\eps>0$; by ``shrinking $I$'', we mean taking $\eps'\in(0,\eps]$ and redefining $I:=(-\eps',\eps')$.

Throughout this work, we will check whether a smooth map is a submersion by looking at properties of curves; in particular, whether they admit local lifts.  This allows us to avoid explicitly dealing with tangent bundles.

\begin{definition}\label{D:LCL}
    A function $f\colon M \to N$ satisfies \textbf{the local curve lifting (LCL) condition} if for any curve $p\colon I\to N$ and $x\in M$ satisfying $f(x)=p(0)$, after possibly shrinking $I$, there exists a (smooth) lift $q\colon I\to M$ of $p$ (restricted to the redefined $I$) through $x$ with respect to $f$.   Explicitly, $q$ satisfies $f\circ q=p$ and $q(0)=x$.
\end{definition}

The LCL condition allows us to identify submersions and diffeomorphisms as follows.  
	\begin{lemma}\label{l:subd}
		Let $M$ and $N$ be manifolds. 
		\begin{enumerate}
		
			\item\label{i:local subd} A smooth surjection $f\colon M\to N$ is a surjective submersion if and only if it satisfies the LCL condition.
			
			\item\label{i:diffeom} A smooth map $f\colon M\to N$ is a diffeomorphism if and only if $f$ is bijective and satisfies the LCL condition.  Moreover, the LCL condition can be relaxed in this case to finding a lift through \emph{any} point.
		\end{enumerate}
	\end{lemma}

	 We will also  make use of the following well-known fact about fibred products of manifolds, which follows from the Transversality Theorem \cite[Theorem 6.30]{lee}. 

	\begin{lemma}\label{l:fibred product}
		Let $f\colon M\to P$ and $g\colon N\to P$ be smooth maps between manifolds, in which $f$ is a surjective submersion.  Then $M\ftimes{f}{g}N$ is a manifold and $\pr_2$ is a surjective submersion.
	\end{lemma}


With this set, we consider our  basic context of  Lie groupoids.
	\begin{definition}\label{d:lie gpds} \cite{moerdijk}
	The $2$-category of \textbf{ Lie groupoids}, denoted   $\LieGpoid$, has objects which are groupoids $\calG = (\calG_1\toto\calG_0)$ in which all structure maps are smooth:   source and target maps $s, t\colon  \calG_1 \to \calG_0$ (which are additionally required to be submersive); the unit map $u\colon  \calG_0 \to \calG_1$;  and the  inverse map $\inv\colon  \calG_1 \to \calG_1$ (where we indicate the inverse $\inv(g)$ by $g^{-1}$).  The $1$-cells are smooth functors; that is the maps on both arrows and objects are smooth.   The $2$-cells are natural transformations between smooth functors that are defined by a smooth map.  
	
	\end{definition}
For the remainder of this paper, unless stated otherwise, we assume that all of our groupoids are Lie groupoids and that all functors and natural transformations are smooth.  \\  

In the study of Lie groupoids, especially orbifolds, actions of groupoids, and stacks, the notion of Morita equivalence is paramount.  We recall its definition:

	\begin{definition} \label{d:ess equiv}  \cite{moerdijk}
		A functor $\varphi\colon\calG\to\calH$ is a \textbf{weak equivalence} (sometimes  called an \textbf{equivalence} or an \textbf{essential equivalence} in the literature) if it satisfies the following two conditions: 
		\begin{enumerate}
			\item\label{def:i:ess surj} \textbf{Smooth Essential Surjectivity:} The induced  map $$\ES{\varphi}\colon\calG_0\ftimes{\varphi_0}{t}\calH_1\to\calH_0\colon (x,h)\mapsto s(h)$$ is a surjective submersion.
			\item\label{def:i:ff} \textbf{Smooth Fully Faithfulness:} The induced map $$\Ff\varphi\colon\calG_1\to\calG_0^2\ftimes{\varphi_0^2}{(s,t)}\calH_1\colon g\mapsto (s(g),t(g),\varphi(g))$$ is a diffeomorphism.
		\end{enumerate}
		We  will denote a weak equivalence with the symbol $\overset{\simeq}{\longrightarrow}$, and denote the class of all weak equivalences in $\LieGpoid$ by $\we$.  
        We say that two Lie groupoids $\calG$ and $\calH$ are \textbf{Morita equivalent} if there is a generalised morphism $\calG\underset{\varphi}{\ot}\calK\underset{\psi}{~\overset{\simeq}{\longrightarrow}~}\calH$ in which {\em both} $\varphi$ and $\psi$ are weak equivalences.
	\end{definition}

We have the following properties of weak equivalences.

    \begin{lemma}\label{l:weak equiv} Let $\varphi\colon\calG\to\calH$ and $\varphi'\colon\calH\to\calK$ be functors in $\LieGpoid$.
        \begin{enumerate}
            \item\label{i:3for2}  If any two of $\varphi$, $\varphi'$, and $\varphi'\circ\varphi$ are weak equivalences, then so is the third.	

            \item\label{i:ssw} If $\varphi$ is smoothly fully faithful and $\varphi_0$ is a surjective submersion, then $\varphi$ is a weak equivalence.

            \item\label{i:ff} A functor $\varphi\colon\calG\to\calH$ in $\LieGpoid$ is smoothly fully faithful if and only if for any functors $\psi,\psi'\colon\calK\to\calG$ and natural transformation $\eta\colon\varphi\circ\psi\Rightarrow\varphi\circ\psi'$, there exists a unique natural transformation $\eta'\colon\psi\Rightarrow\psi'$ such that $\eta=\varphi\eta'$.
        \end{enumerate}
    \end{lemma}

        In fact, the weak equivalences as in Item~\ref{i:ssw} are important enough to this paper to warrant a definition.
        
    \begin{definition}\label{d:ssw}
        A weak equivalence that is a surjective submersion on objects is called a \textbf{surjective submersive weak equivalence}, often denoted using $\overset{\simeq}{\longtwoheadrightarrow}$.  (They are also called \textbf{Morita fibrations} in \cite{dHF}.)  We denote the class of surjective submersive weak equivalences by $\ssw$.
    \end{definition}

    \begin{proof}[Proof of Lemma~\ref{l:weak equiv}:]
        Item \ref{i:3for2} is the so-called 3-for-2 property, as proved in \cite[Lemma 8.1]{pronk-scull-cor}.

        Item~\ref{i:ssw} appears in \cite[Subsection 6.1]{dHF}; alternatively, the reader can check that $\ES\varphi$ satisfies the LCL condition and then apply Lemma~\ref{l:subd}.

        For Item \ref{i:ff}, suppose $\varphi$ is smoothly fully faithful.  Fix functors $\psi,\psi'\colon\calK\to\calG$ and natural transformation $\eta\colon\varphi\circ\psi\Rightarrow\varphi\circ\psi'$.  Define $\eta'\colon\calK_0\to\calG_1$ by $\eta'(z):=\Ff\varphi^{-1}(\psi(z),\psi'(z),\eta(z))$; this is well-defined and smooth.  The fact that  $\varphi\eta'=\eta$ follows from the construction, and uniqueness follows from smooth fully faithfulness of $\varphi$.
			
		Conversely, suppose for any functors $\psi,\psi'\colon\calK\to\calG$ and natural transformation $\eta\colon\varphi\circ\psi\Rightarrow\varphi\circ\psi'$, there exists a unique natural transformation $\eta'\colon\psi\Rightarrow\psi'$ such that $\eta=\varphi\eta'$. Let $\calK$ be the trivial Lie groupoid of a point $*\toto*$, and fix $(x,x',h)\in\calG_0^2\ftimes{\varphi^2}{(s,t)}\calH_1$.  Set $\psi\colon\calK\to\calG$ to be the functor sending the point to $x$, and set $\psi'\colon\calK\to\calG$ the functor sending the point to $x'$.  The natural transformation $\eta\colon\varphi\circ\psi\Rightarrow\varphi\circ\psi'$ sending the point to $h$ factors uniquely as $\varphi\eta'$.  Thus, $\Ff{\varphi}(\eta'(*))=(x_1,x_2,h)$, and $\Ff{\varphi}$ is bijective.  
				
		Now we show that $\Ff\varphi$ satisfies the LCL condition.  Fix a curve $p=(x_\tau, x'_\tau, h_\tau)\colon I\to\calG_0^2\ftimes{\varphi^2}{(s,t)}\calH_1$.  Let $\calK$ to be the trivial Lie groupoid $I\toto I$; let $\psi,\psi'\colon\calK\to\calG$ be the functors defined by $\psi_0= x_\tau$ and $\psi'_0 = x_\tau'$; and let $\eta\colon  \varphi\circ\psi \Rightarrow \varphi\circ\psi'$ defined by the natural transformation $h_\tau$.    There is a unique $\eta'\colon\psi\Rightarrow\psi'$ defined by $\eta'(\tau)=g_\tau\in\calG_1$ such that $\varphi\eta'=\eta$, and so $\varphi(g_\tau) = h_\tau$.  Thus $\Ff\varphi(g_\tau)=p$, showing that $g_\tau$ gives the  required lift.    By Item~\ref{i:diffeom} of Lemma~\ref{l:subd}, it follows that $\Ff\varphi$ is a diffeomorphism, establishing smooth fully faithfulness of $\varphi$. This proves Item \ref{i:ff}.        
    \end{proof}


We next consider  pullbacks of Lie groupoids.  
	\begin{definition}\label{d:strict pullback}
		Let $\varphi\colon\calG\to\calK$ and $\psi\colon\calH\to\calK$ be functors.  The \textbf{strict pullback} of $\varphi$ and $\psi$ is the groupoid $\calG\ftimes{\varphi}{\psi}\calH$, whose object and arrow spaces are the corresponding fibred products of the object and arrow spaces of $\calG$ and $\calH$, respectively,  with two projection functors $\pr_1$ and $\pr_2$.   
	\end{definition}

The strict pullback may not be a Lie groupoid in general.  The following proposition provides a sufficient condition for when it is a Lie groupoid.  

	\begin{proposition}\label{p:strict pullback}
		Let $\varphi\colon\calG\to\calK$ and $\psi\colon\calH\to\calK$ be functors in which $\varphi\in\ssw$.  Then $\calG\ftimes{\varphi}{\psi}\calH$ is a Lie groupoid and $\pr_2\in\ssw$.
	\end{proposition}

	\begin{proof}
		By Lemma \ref{l:fibred product}, the object and arrow spaces of $\calG\ftimes{\varphi}{\psi}\calH$ are manifolds. We verify that the source map of the pullback groupoid is a surjective submersion using the LCL condition.  Fix a curve $p=(x_\tau, y_\tau)\colon I\to\calG_0\ftimes{\varphi_0}{\psi_0}\calH_0$ and let $(g_0,h_0)\in\calG_1\ftimes{\varphi_1}{\psi_1}\calH_1$ such that $s(g_0,h_0)=(x_0, y_0)$.  After shrinking $I$, there is a lift $h_\tau\colon I\to\calH_1$ of $y_\tau$ through $h_0$ such that $s(h_\tau) = y_\tau$.  Then $t(\psi_1(h_\tau))$ defines a curve $I \to \calK_0$, and since $\varphi_0$ is a surjective submersion, after shrinking $I$ again there is a lift $x'_\tau\colon I\to\calG_0$ of this curve  through $t(g_0)$ with $\varphi(x'_\tau) = t(\psi(h_\tau))$. The curve $g_\tau:=\Ff{\varphi}^{-1}(x_\tau, x'_\tau,\psi_1(h_\tau))\colon I \to\calG_1$ is a lift of $x_\tau$ through $g_0$ such that $s(g_\tau)=x_\tau$.  Thus $(g_\tau, h_\tau)\colon I\to\calG_1\ftimes{\varphi_1}{\psi_1}\calH_1$ is well-defined, and is the desired lift of $(x_\tau, y_\tau)$ through $(g_0,h_0)$ such that $s(g_\tau, h_\tau) = (x_\tau, y_\tau)$  verifying the LCL condition for the source map of the pullback groupoid.  By Item~\ref{i:local subd} of Lemma~\ref{l:subd}, the source map is a surjective submersion, from which it follows that the target map is as well.  Thus $\calG\ftimes{\varphi}{\psi}\calH$ is a Lie groupoid.

		Next we show that $\pr_2$ is in $\ssw$.  By Lemma~\ref{l:fibred product}, the map $(\pr_2)_0\colon(\calG\ftimes{\varphi}{\psi}\calH)_0\to\calH_0$ is a surjective submersion.    Since $\Ff{\varphi}$ is a diffeomorphism, it follows that $\Ff{\pr_2}$ is bijective.   Let $$p=((x_\tau, y_\tau), (x_\tau', y_\tau'),h_\tau)\colon I\to(\calG\ftimes{\varphi}{\psi}\calH)_0^2\ftimes{\varphi_0^2}{(s,t)}\calH_1$$ be a curve. Define the curve $g_\tau:=\Ff{\varphi}^{-1}((x_\tau, x'_\tau),\psi_1(h_\tau))$.  Then $(g_\tau,h_\tau)$ defines the desired lift of the curve $p$ and $\Ff{\pr_2}$ satisfies the   LCL condition.  Hence $\Ff{\pr_2}$ is a surjective submersion.  By Item~\ref{i:diffeom} of Lemma~\ref{l:subd}, $\Ff{\pr_2}$ is a diffeomorphism.  By Item~\ref{i:ff} of Lemma~\ref{l:weak equiv}, $\pr_2$ is in $\ssw$.
	\end{proof}
	
We will also be using the weak pullback of groupoids.  

	\begin{definition}\label{d:weak pullback}
	 \cite{moerdijk} Let $\varphi\colon\calG\to\calK$ and $\psi\colon\calH\to\calK$ be functors in $\LieGpoid$. The \textbf{weak pullback} of $\varphi$ and $\psi$ is the groupoid $		\calG\wtimes{\varphi}{\psi}\calH$, whose object space is
$$(\calG\wtimes{\varphi}{\psi}\calH)_0:=\left\{(x,k,y)\in\calG_0\times\calK_1\times\calH_0\mid \varphi(x)\overset{k}{\curvearrowright}\psi(y)\right\},$$
and arrow space is
$$(\calG\wtimes{\varphi}{\psi}\calH)_1:=\left\{(g,k,h)\in\calG_1\times\calK_1\times\calH_1\mid \varphi(s(g))\overset{k}{\curvearrowright}\psi(s(h))\right\}.$$

The weak pullback comes equipped with two projection functors $\pr_1$ and $\pr_3$ to $\calG$ and $\calH$, respectively, and the natural transformation  $\PR_2\colon\varphi\circ\pr_1\Rightarrow\psi\circ\pr_3$.
	\end{definition}

In general, the weak pullback may not be a Lie groupoid.  The following proposition, which is \cite[Proposition 5.12(iv)]{MM}, gives a sufficient condition for when it is a Lie groupoid.

	\begin{proposition}\label{p:weak pullback}
		Let $\varphi\colon\calG\to\calK$ and $\psi\colon\calH\to\calK$ be functors in which $\varphi$ is a weak equivalence.  Then $\calG\wtimes{\varphi}{\psi}\calH$ is a Lie groupoid and $\pr_3\in\ssw$.
	\end{proposition}

The following lemma shows how surjective submersive weak equivalences interact with natural transformations.   This property is called ``co-fully faithfulness'' by Roberts \cite[Definition 2.12]{roberts2021} and Pronk-Scull \cite[Definition 5.1]{pronk-scull:bicategories}.

	\begin{lemma}\label{l:co-ff}
		Given $\varphi\colon\calG\to\calH$ be in $\ssw$, for any functors $\psi,\psi'\colon\calH\to\calK$ and natural transformation $\eta\colon\psi\circ\varphi\Rightarrow\psi'\circ\varphi$, there exists a unique natural transformation $\eta'\colon\psi\Rightarrow\psi'$ such that $\eta=\eta'\varphi$.
	\end{lemma}

	\begin{proof}
		Fix functors $\psi,\psi'\colon\calH\to\calK$ and natural transformation $\eta\colon\psi\circ\varphi\Rightarrow\psi'\circ\varphi$.  Define $\eta'\colon\calH_0\to\calK_1$ by $\eta'(y):=\eta(x)$, where $x\in\varphi_0^{-1}(y)$: since $\varphi_0$ is surjective,  $\varphi_0^{-1}(y)$ is non-empty.  Suppose $\varphi_0(x_1)=\varphi_0(x_2)$.  Since $\varphi$ is a weak equivalence, there exists an arrow $g=\Ff\varphi^{-1}(x_1,x_2,u_{\varphi(x_1)})$ from $x_1$ to $x_2$.  Naturality gives the following commutative diagram
		$$\xymatrix{
			\psi\circ\varphi(x_1) \ar[r]^{\eta_{x_1}} \ar[d]_{\psi\circ\varphi(g)} & \psi'\circ\varphi(x_1) \ar[d]^{\psi'\circ\varphi(g)} \\
			\psi\circ\varphi(x_2) \ar[r]_{\eta_{x_2}} & \psi'\circ\varphi(x_2).
		}$$
Since $\varphi(g)=u_{\varphi(x_1)}$, we have $\eta_{x_1}=\eta_{x_2}$, and so $\eta'$ is well-defined.    By construction $\eta=\eta'\varphi$.

		To show that $\eta'$ is smooth, fix a curve $p=y_\tau\colon I\to\calH_0$.  After shrinking $I$, there exists a curve $x_\tau\colon I\to\calG_0$ such that  $y_\tau=\varphi(x_\tau)$, since $\varphi$ is a surjective submersion.  Since $\eta'(y_\tau) = \eta(x_\tau)$, we conclude that $\eta'(y_\tau)$ is a curve in $\calK_1$.   By Boman's Lemma \cite{boman}, $\eta'$ is smooth.  The naturality of $\eta'$ follows from  that  of $\eta$.  Finally, uniqueness follows from the construction.
	\end{proof}
The following identifies weak equivalences using a property  called ``$J$-locally split'' in \cite[Definition 3.22]{roberts2021}; in our case, $J=\ssw$.
	\begin{lemma}\label{l:loc split}
		A functor $\varphi\colon\calG\to\calH$ is a weak equivalence if and only if it is smoothly fully faithful and $\ssw$-locally split: \emph{i.e.}\ there exist a functor $\psi\colon \calK\to \calH$ in $\ssw$, a functor $\sigma\colon \calK\to \calG$, and a natural transformation $\eta\colon \varphi\circ\sigma\Rightarrow\psi$.   
	\end{lemma}
 
	\begin{proof}
		Suppose $\varphi$ is a weak equivalence.  Choose $\calK:=\calG\wtimes{\varphi}{\id_\calH}\calH$ and $\psi=\pr_3$, $\sigma=\pr_1$ and $\eta=\PR_2$.  Then $\calK$ is a Lie groupoid and $\psi \in \ssw$ by Proposition~\ref{p:weak pullback}.
			
	Conversely, suppose $\varphi$ is smoothly fully faithful, and there exist $\psi\colon \calK\to \calH$ in $\ssw$, a functor $\sigma\colon \calK\to \calG$, and a natural transformation $\eta\colon \varphi\circ\sigma\Rightarrow\psi$.  We will verify that $\ES\varphi$ is a surjective submersion.

		Let $y\in\calH_0$.  Since $\psi_0$ is surjective, there exists $z\in\calK_0$ such that $\psi(z)=y$.  Then $(\sigma(z),\eta(z)^{-1})\in\calG_0\ftimes{\varphi}{t}\calH_1$ and $\ES{\varphi}(\sigma(z),\eta(z)^{-1})=y$.  Thus $\ES{\varphi}$ is surjective.
		To show $\ES\varphi$ is a surjective submersion, fix a curve $p = y_\tau\colon I\to\calH_0$. Let $(x_0,h_0)\in\calG_0\ftimes{\varphi}{t}\calH_1$ such that $\ES\varphi(x_0,h_0)=y_0$.    Since $\psi_0$ is a surjective submersion, after shrinking $I$, there is a curve $z_\tau\colon I\to\calK_0$ such that $\psi(z_\tau) =y_\tau$.  Since $\varphi$ is smoothly fully faithful, $\Ff\varphi$ is a diffeomorphism and so  we define  $g_0 \in\calG_1$  by $\Ff{\varphi}^{-1} (\sigma(z_0),x_0,h_0\cdot\eta(z_0)).$   Since source maps of Lie groupoids are surjective submersions, we can lift the curve $\sigma (z_\tau)$   after shrinking $I$ to get a curve $g_\tau\colon I\to\calG_1$ through $g_0$ such that $s(g_\tau) =\sigma(z_\tau)$.  Then the desired lift of the  curve $p$ through $(x_0, h_0)$ is defined by  $(t(g_\tau),\varphi(g_\tau)\cdot\eta(z_\tau)^{-1})$.    Thus  $\ES\varphi$ satisfies the LCL condition and is a surjective submersion, completing the verification that  $\varphi$ is a weak equivalence.
	\end{proof}

\section{Localising Lie Groupoids at Weak Equivalences}\label{s:localisation}

In this section we recall how to construct  a localised bicategory which inverts weak equivalences in the $2$-category $\LieGpoid$, giving  us a formal mechanism for working with Morita equivalence classes of groupoids.    We will use a bicategory of fractions construction in which the objects are still the Lie groupoids of $\LieGpoid$, but the 1- and 2-cells are adjusted.  In particular, the arrows of the bicategory of fractions will be given by so-called ``spans'' of arrows of $\LieGpoid$, so that we add inverse arrows for any weak equivalences and make all weak equivalences  (and hence all Morita equivalences) into isomorphisms of the localised bicategory.

We outline  two constructions of this localisation: the first is the bicategory of fractions defined  by \cite{pronk}, and a second related but smaller construction based on so-called anafunctors by \cite{roberts2021}.  

We begin by recalling the construction and properties of the localized bicategory of \cite{pronk}.

 \begin{definition}\label{d:gen morph}  The localized bicategory $\LieGpoid[\we^{-1}]$ is defined by:
 \begin{itemize}
     \item The objects of $\LieGpoid[\we^{-1}]$ are the same as the objects of  $\LieGpoid$.
 
	\item The arrows are defined by spans of arrows: 	a \textbf{generalised morphism} between Lie groupoids $\calG$ and $\calH$ is a Lie groupoid $\calK$ and two  functors  $\calG\underset{\varphi}{\ot}\calK\underset{\psi}{\longrightarrow}\calH$ in which $\varphi$ is a weak equivalence .
\item The \textbf{identity generalised morphism} of  $\calG$ in $\LieGpoid[\we^{-1}]$ is given by $\calG\overset{=}{\longleftarrow}\calG\overset{=}{\longrightarrow}\calG$.
\item 	Let $\calG\underset{\varphi}{\ot}\calK\underset{\psi}{\longrightarrow}\calH$ and $\calH\underset{\chi}{\ot}\calL\underset{\omega}{\longrightarrow}\calI$ be generalised morphisms.    Define their \textbf{composition} to be the generalised morphism $$\calG\underset{\varphi\circ\pr_1}{\ot}\calK\wtimes{\psi}{\chi}\calL\underset{\omega\circ\pr_3}{\longrightarrow}\calI.$$
\item  2-cells between generalised morphisms are  detailed below in Definition \ref{d:equiv gen morph}.
\end{itemize}
	\end{definition}
	
It follows from Proposition~\ref{p:weak pullback} and Item~\ref{i:3for2} of Lemma~\ref{l:weak equiv} that the composition of two generalised morphisms is a generalised morphism.
	The composition of generalised morphisms is an associative operation up to a canonical isomorphism.
	
We think of a generalised morphism $\calG\ot\calK\to\calH$ as replacing $\calG$ with a weakly equivalent Lie groupoid $\calK$ which admits the left functor $\calK\to\calH$. Thus we consider weakly equivalent groupoids to be interchangeable. 
Weakly equivalent groupoids are always Morita equivalent, as there is a generalised morphism between them using an identity morphism as one leg. 
A Morita equivalence $\calG\underset{\varphi}{\ot}\calK\underset{\psi}{~\overset{\simeq}{\longrightarrow}~}\calH$ is invertible in  $\LieGpoid[\we^{-1}]$, with  inverse defined by $\calH\underset{\psi}{\ot}\calK\underset{\varphi}{~\overset{\simeq}{\longrightarrow}~}\calG$.

There may be many different choices of  groupoids weakly equivalent to  $\calG$, and we want to recognise when two choices of generalised morphism carry the same geometric information.   Thus we define the following 2-cells.   

	\begin{definition}
    \label{d:equiv gen morph}
		Given two generalised morphisms $\calG\underset{\varphi}{\ot}\calK\underset{\psi}{\longrightarrow}\calH$ and $\calG\underset{\varphi'}{\ot}\calK'\underset{\psi'}{\longrightarrow}\calH$ we consider a generalised morphism of the form  $\calK\underset{\alpha}{\ot}\calL\underset{\alpha'}{\overset{\simeq}{\longrightarrow}}\calK'$ in which both functors are weak equivalences, along with two natural transformations $\eta_1\colon\varphi\circ\alpha\Rightarrow\varphi'\circ\alpha'$ and $\eta_2\colon\psi\circ\alpha\Rightarrow\psi'\circ\alpha'$ making the following diagram $2$-commute:
		\begin{equation}\label{e:equiv gen morph 2cell}
		\xymatrix{& & \calK \ar[dll]_{\varphi}^{\simeq} \ar[drr]^{\psi} & & \\
		\calG & \Downarrow~\eta_1 & \calL \ar[u]^{\alpha}_{\simeq} \ar[d]_{\alpha'}^{\simeq} & \Downarrow~\eta_2 & \calH \\
		 & & \calK'. \ar[ull]^{\varphi'}_{\simeq} \ar[urr]_{\psi'} & & \\
		}
		\end{equation}
		We denote a diagram \eqref{e:equiv gen morph 2cell} by the quadruple $(\alpha,\alpha',\eta_1,\eta_2)$. A \textbf{$2$-cell} from $\calG\underset{\varphi}{\ot}\calK\underset{\psi}{\longrightarrow}\calH$ to $\calG\underset{\varphi'}{\ot}\calK'\underset{\psi'}{\longrightarrow}\calH$ is an equivalence class of such diagrams of form  \eqref{e:equiv gen morph 2cell}:     $(\alpha,\alpha',\eta_1,\eta_2)$ is equivalent to another such diagram $(\beta,\beta',\mu_1,\mu_2)$ if there exists a third  generalised morphism $\calL\underset{\gamma}{\ot}\calN\underset{\gamma'}{\overset{\simeq}{\longrightarrow}}\calM$ and natural transformations $\nu_1\colon\alpha\circ\gamma\Rightarrow\beta\circ\gamma'$ and $\nu_2\colon\alpha'\circ\gamma\Rightarrow\beta'\circ\gamma'$ such that 
  
		\begin{equation}\label{e:equiv gen morph}
(\mu_1\gamma')\circ(\varphi\nu_1)=(\varphi'\nu_2)\circ(\eta_1\gamma) \quad \text{and} \quad (\mu_2\gamma')\circ(\psi\nu_1)=(\psi'\nu_2)\circ(\eta_2\gamma).
		\end{equation}   
  
We denote the 2-cell given by the equivalence class by   $[\alpha,\alpha',\eta_1,\eta_2]$.
		The \textbf{identity $2$-cell} of a generalised morphism $\calG\ot\calK\to\calH$ is given by $[\id_\calK,\id_\calK,\ID_\varphi,\ID_\psi]$ where $\ID$ represents the identity natural transformation.   
	\end{definition}

		Vertical compositon of 2-cells in $\LieGpoid[\we^{-1}]$ is performed using a weak pullback, and horizontal compositions are defined using whiskering operations.  Unitors for this bicategory are defined using projection maps.  
Explicit descriptions of 2-cell compositions in $\LieGpoid[\we^{-1}]$
		can be found in  \cite[Subsection 2.3]{pronk} or \cite[Section 3]{pronk-scull:bicategories}.

The  bicategory $\LieGpoid[\we^{-1}]$ inverts all weak equivalences and satisfies the universal property of a localisation.

	\begin{proposition}\label{p:liegpoid bicateg of frac}  \cite[Section 2]{pronk} The bicategory  $\LieGpoid[\we^{-1}]$ satisfies the universal property of a localisation:  any functor from $\LieGpoid$ to another bicategory which takes $1$-cells in $\we$ to invertible $1$-cells will factor through this localised bicategory  $\LieGpoid[\we^{-1}]$.  
	\end{proposition}



On a practical level,  this bicategory can be hard to work with on the $2$-cell level, since the $2$-cells are defined as equivalence classes of diagrams.   We now describe an alternate localised bicategory $\ana$ developed in \cite{roberts2012,roberts2016,roberts2021} which is ``smaller'' than  $\LieGpoid[\we^{-1}]$ described in the previous section.    The fact that this smaller construction applies in the category of Lie groupoids is well-known to experts but tracking down exact references has proved difficult.

The starting point to understanding the difference between $\LieGpoid[\we^{-1}]$ and $\ana$ is in looking at the  subclass $\ssw$ of $\we$ comprising surjective submersive weak equivalences.  \emph{A priori} we will create a localisation which inverts only this subclass $\ssw$.  It will turn out that the resulting localised bicategories are equivalent.

We start creating $\ana$ using objects of $\LieGpoid$ as before, but look only at generalised morphisms which use a surjective submersive weak equivalence as their left leg.  These are called anafunctors in 
\cite{roberts2021}. 

	\begin{definition}\label{d:ana gen morph}  The category $\ana$ is defined as follows:
 \begin{itemize}
     \item objects  are the same as the objects of  $\LieGpoid$.
     \item a $1$-cell  is given by an  \textbf{anafunctor}, a generalised morphism $\calG\underset{\varphi}{\ot}\calK\underset{\psi}{\longrightarrow}\calH$  such that  $\varphi \in \ssw$.
     \item  The identity generalised morphism of Definition \ref{d:gen morph} is an anafunctor, and so it defines the  \textbf{identity anafunctor} of $\calG$.
     \item Compositon is defined using the strict pullback:  Let $\calG\underset{\varphi}{\ot}\calK\underset{\psi}{\longrightarrow}\calH$ and $\calH\underset{\chi}{\ot}\calL\underset{\omega}{\longrightarrow}\calI$ be anafunctors.  Define their \textbf{composition} to be the anafunctor $$\calG\underset{\varphi\circ\pr_1}{\ot}\calK\ftimes{\psi}{\chi}\calL\underset{\omega\circ\pr_3}{\longrightarrow}\calI.$$
     \item the 2-cells are defined by natural transformations (not equivalence classes of diagrams):   Given two  anafunctors $\calG\underset{\varphi}{\ot}\calK\underset{\psi}{\longrightarrow}\calH$ and $\calG\underset{\varphi'}{\ot}\calK'\underset{\psi'}{\longrightarrow}\calH$   a \textbf{$2$-cell} between them is a natural transformation $\eta$ making the following diagram $2$-commute.  
		$$\xymatrix{
		 & \calK \ar[dr]^{\psi} & \\
		\calK\ftimes{\varphi}{\varphi'}\calK' \ar@{>>}[ur]^{\pr_1}_{\simeq} \ar@{>>}[dr]_{\pr_2}^{\simeq} & \overset{\eta}{\Rightarrow} & \calH \\
		 & \calK' \ar[ur]_{\psi'} & \\
		}$$
 \item The \textbf{identity $2$-cell} of an anafunctor $\calG\ot\calK\to\calH$ is given by the natural transformation $$\iota_{\calG\leftarrow\calK\to\calH}\colon(\calK\ftimes{\varphi}{\varphi}\calK)_0\to\calH_1\colon(y_1,y_2)\mapsto\psi(\Ff{\varphi}^{-1}((y_1,y_2),u_{\varphi(y_1)})).$$
 \end{itemize}
		
	\end{definition}

It follows from Proposition~\ref{p:strict pullback}, Item~\ref{i:3for2} of Lemma~\ref{l:weak equiv}, and the fact that the composition of surjective submersions is again a surjective submersion that the composition of two anafunctors is an anafunctor.  Similar to the case of generalised morphisms, composition of anafunctors is an associative operation, again up to a canonical isomorphism.    
	
We have chosen a canonical representative of a $2$-cell defined by a particular natural transformation which will  \emph{be} the $2$-cell in our new bicategory $\ana$. Thus the  $2$-cells are  actual natural transformations rather than equivalence classes as in $\LieGpoid[\we^{-1}]$; these  anafunctor 2-cells can be drawn in the diagram form \eqref{e:equiv gen morph 2cell} of the representatives of 2-cells of $\LieGpoid[\we^{-1}]$ of Definition \ref{d:equiv gen morph}, by adding the left side of the diagram with the trivial natural transformation.   However, in $\ana$, the vertical and horizontal compositions are {\it not} the usual composition of natural transformations.   We do not need the details of these compositions until  Section \ref{s:eqvt-to-full}, so we will defer the relevant details.


Both of the constructions we have just outlined create localisations of $\LieGpoid$ at $\we$.  To compare $\ana$ to $\LieGpoid[\we^{-1}]$, we again apply a result of Roberts; be warned that what is called a ``weak equivalence'' in his paper \cite{roberts2021} is defined there to be a  functor that is $\ssw$-locally split and representably fully faithful.  However, by Item~\ref{i:ff} of Lemma~\ref{l:weak equiv} and Lemma~\ref{l:loc split}, this is equivalent to smooth essential surjectivity and smooth fully faithfulness, and so our notion of weak equivalence coincides with his.  Thus we have  \cite[Theorem 3.24]{roberts2021}:

	\begin{proposition}\label{p:ana bicat2}
		The inclusion $\ana \to \LieGpoid[\we^{-1}]$ is an equivalence of bicategories, where this inclusion takes a $2$-cell to its equivalence class.   
	\end{proposition}

 Proposition \ref{p:ana bicat2} implies that any generalised morphism $\calG\underset{\varphi}{\ot}\calK\underset{\psi}{\longrightarrow}\calH$ admits a $2$-cell from itself to an anafunctor; in the proof, this is the anafunctor $\calG\underset{\pr_1}{\ot}\calG\wtimes{\id_\calG}{\varphi}\calK\underset{\psi\circ\pr_3}{\longrightarrow}\calH$.  
 It also implies that there is a $2$-cell from a composition of generalised morphisms in $\LieGpoid[\we^{-1}]$, defined by the weak pullback, to the corresponding composition in $\ana$ using the strict pullback.  This can be constructed explicitly as a vertical composition using the following, which we use in Section~\ref{s:eqvt-to-full}.

	\begin{proposition}\label{p:comp gen morph}  Let $\calG\underset{\varphi}{\ot}\calK\underset{\psi}{\longrightarrow}\calH$ and $\calH\underset{\chi}{\ot}\calL\underset{\omega}{\longrightarrow}\calI$ be generalised morphisms where $\chi\in \ssw$ so that the second generalised morphism is  an anafunctor.  
		  There is a $2$-cell from the composition defined using the weak pullback in $\LieGpoid[\we^{-1}]$ of Definition \ref{d:gen morph}  to the composition defined using the strict pullback in $\ana$ defined in Definition \ref{d:ana gen morph}.
	\end{proposition}
	
		  \begin{proof}
		 We  define a  $2$-cell between the two compositions using the  following  strictly commutative diagram:
		$$\xymatrix{
		 & & \calK\wtimes{\psi}{\chi}\calL \ar[dll]_{\varphi\circ\pr_1}^{\simeq} \ar[drr]^{\omega\circ\pr_3} & & \\
		\calG & \circlearrowright & \calK\ftimes{\psi}{\chi}\calL \ar[u]^{\ic}_{\simeq} \ar@{=}[d] & \circlearrowright & \calI \\
		 & & \calK\ftimes{\psi}{\chi}\calL. \ar[ull]^{\varphi\circ\pr_1}_{\simeq} \ar[urr]_{\omega\circ\pr_2} & & \\
		}$$
Here, $\ic$ is the inclusion functor defined on objects by $\ic_0\colon(x,y)\mapsto(x,u_{\psi(x)},y)$ and on arrows by $\ic_1\colon(k,\ell)\mapsto(k,u_{s\circ\psi(k)},\ell)$. We check that $\ic$ is  a weak equivalence.   It is straightforward to check that $\Ff{\ic}$ is bijective.  To show it is a diffeomorphism, we check the LCL condition:  a curve  $p\colon I\to(\calK\ftimes{\psi}{\chi}\calL)_0^2\ftimes{\ic^2}{(s,t)}(\calK\wtimes{\psi}{\chi}\calL)_1$ has the form $p(t) = ((x_\tau, y_\tau), (x'_\tau, y'_\tau), (k_\tau, u_{x_\tau}, \ell_\tau))$ and so   $p=\Ff{\ic}(k_\tau, \ell_\tau) $ and $(k_\tau, \ell_\tau)$ gives the desired lift.     By Item~\ref{i:diffeom} of Lemma~\ref{l:subd}, $\Ff{\ic}$ is a diffeomorphism.

		Let $(x,h,y)\in(\calK\wtimes{\psi}{\chi}\calL)_0$.  Since $\chi_0$ is surjective, there exists $y'\in\calL_0$ such that $\chi(y')=\psi(x)$.  Define $\ell:=\Ff{\chi}^{-1}(y,y',h^{-1})$.  Then $\ES{\ic}((x,y'),(u_{x},h,\ell))=(x,h,y)$, and so $\ES{\ic}$ is surjective.

		To show that $\ES{\ic}$ is a surjective submersion,  we again use the LCL condition:   let $p=(x_\tau, h_\tau, y_\tau)\colon I\to (\calK\wtimes{\psi}{\chi}\calL)_0 $ be a curve  and  suppose  we have a point $((x_0',y_0'), (k_0, h_0, \ell_0)) \in ((\calK \ftimes{\psi}{\chi} \calL)_0)\ftimes{\ic}{t}  (\calK\wtimes{\psi}{\chi}\calL)_1$ sent by $\ES{inc}$ to $(x_0,h_0,y_0)$.  By definition, we must have that $t(k_0) = x_0'$, $t(\ell_0) = y_0'$, and $\chi(\ell_0) h_0 \psi(k_0)^{-1} = u_{\chi(y'_0)}$.
		Since $\calK$ is a Lie groupoid, its source map is a surjective submersion and so after shrinking $I$, there is a curve $k_\tau\colon I\to\calK_1$ through $k_0$ with  $s(k_\tau) = x_\tau$;  denote the target $t(k_\tau) = x'_\tau$, and note that $\psi(x'_0) = t(\psi(k_0)) =  t(\chi(\ell_0)h_0)  = t \chi(\ell_0)$.   Next, since  $\chi_0$ is a surjective submersion we can lift  $t (\psi (k_\tau)) = \psi(x_\tau') $ through $t(\ell_0)$ to get $y_\tau'$.   Thus we have $(y_\tau, y'_\tau, \psi(k_\tau)h_\tau^{-1}) \in \calL_0^2 \ftimes{\chi^2}{(s,t)} \calH_1$ and so we can define $\ell_\tau = \Ff{\chi}^{-1}(y_\tau, y'_\tau, \psi(k_\tau)h_\tau^{-1})$.  Then $q=(x'_\tau, y'_\tau,(k_\tau, h_\tau, \ell_\tau))$ is the desired lift of $p$ with   $\ES{\ic}(q) = p$.

		By Item~\ref{i:3for2} of Lemma~\ref{l:weak equiv}, since $\varphi$ is a weak equivalence, $\varphi\circ\pr_1$ is a weak equivalence (in both instances it appears in the diagram above) provided $\pr_1$ is.  But by Proposition~\ref{p:weak pullback}, $\pr_1\colon\calK\wtimes{\psi}{\chi}\calL\to\calK$ is a weak equivalence since $\chi$ is, and by Proposition~\ref{p:strict pullback} $\pr_1\colon\calK\ftimes{\psi}{\chi}\calL$ is a weak equivalence (in fact, it is in $\ssw$) since $\chi$ is in $\ssw$.   Thus the diagram above represents an equivalence between the two generalised morphisms.
	\end{proof}

\section{Action Groupoids}\label{s:eqactgpd}

Our main interest in this paper is in Lie groupoids which come from the smooth action of a Lie group on a manifold.  For the rest of the paper, we will focus on these, and so  all group actions are assumed to be Lie group actions henceforth.

Recall that the \textbf{action groupoid} of a Lie group action $G \ltimes X$ is defined to be the Lie groupoid with object space $X$; arrow space $G\times X$; source and target given by the second projection map and the action map, resp.; multiplication $(g_1,g_2x)(g_2,x)=(g_1g_2,x)$; unit at $x$ given by $(1_G,x)$; and inverse $(g,x)^{-1}=(g^{-1},gx)$.
	
We are interested in looking at action groupoids with various special properties which commonly come up in contexts such as the study of orbifolds, symplectic geometry, and bundle theory.  Recall that a group action is \textbf{free} if all stabilisers are trivial, \textbf{locally free} if there is a neighbourhood $U$ of $1_G$ in $G$ such that the restriction of the action to $U$ is free, \textbf{transitive} if for each pair $x,y\in X$ there exists $g\in G$ such that $gx=y$, and \textbf{effective} if for each  $g\not= 1_G \in  G$ there exists $x \in X$ such that $gx \not= x$. We will apply these adjectives to both the action and the corresponding action groupoid.  We will also refer to the action as \textbf{compact} (resp.\ \textbf{discrete}) if the corresponding group is compact (resp.\ discrete).  A Lie groupoid $\calG$ is \textbf{proper} if the map $\calG_1\to\calG_0\times\calG_0\colon g\mapsto (s(g),t(g))$ is a proper map.  In particular, if $\calG$ is an action groupoid, then we say that the corresponding action is \textbf{proper}.  Finally, an \textbf{\'etale groupoid} is a Lie groupoid whose source (and hence target) is a local diffeomorphism, and an \textbf{orbifold groupoid} is a Lie groupoid that is Morita equivalent to a proper \'etale groupoid.

        \begin{remark}\label{r:orbifold gpd}
           There are subtle differences in how ``orbifold groupoid'' is defined in the literature.  Our definition above matches that of Pronk-Scull \cite[Definition 2.7]{pronk-scull}.  However, others refer to only proper \'etale groupoids as ``orbifold groupoids''; see, for instance,  \cite{ALR,HM}.  Further, some authors restrict their attention to effective orbifolds \cite{MP,MM}.  It follows from the Slice Theorem and associated Tube Theorem \cite[Theorems 2.3.3, 2.4.1]{DK} that a proper and locally free group action corresponds to an orbifold groupoid as in the definition above. Conversely, if an action groupoid is Morita equivalent to a proper \'etale groupoid, then it is proper and locally free.  This follows from the fact that weak equivalences, and hence Morita equivalences, preserve stabilisers and properness (see, for instance, \cite[Subsection 2.7]{moerdijk}, \cite[Proposition 5.1.5]{del hoyo}, or \cite[Proposition 2.2]{zung}). 
        \end{remark}
        
We will localise a  specified sub-$2$-category of Lie groupoids whose objects are action groupoids that satisfy a desired set of properties $\calP$; this is Proposition~\ref{p:eqvt bicat 2}.   Our selection of properties comes from those defined above.  Specifically,     $\calP$ is any subset (possibly empty) of the following list of properties (acknowledging that some combinations are redundant; see Remark~\ref{r:orbifold gpd}): 
	\begin{equation}\label{e:P}
	\calP \subseteq	\Bigg\{\text{\parbox{3.5in}{\begin{center}free, locally free, transitive, effective, compact, discrete, proper, is an orbifold groupoid\end{center}}}\Bigg\}.
	\end{equation}

    \begin{remark}\label{r:morita equivalent}
        Several of the properties in $\calP$ are known to be Morita invariants; that is, preserved by weak equivalences, and hence Morita equivalences: free, locally free, transitive, proper, and being an orbifold groupoid.  See, for instance, \cite[Theorem 4.3.1, Proposition 5.1.5]{del hoyo} or \cite[Proposition 2.2]{zung}.  This fact will be important in some proofs below.  In fact, we could add the property of being Morita equivalent to any class of groupoids into $\calP$, and all of the following results will continue to hold.
    \end{remark}
 
Similar results to Proposition~\ref{p:eqvt bicat 2} already appear in the literature for a few specific sub-classes of Lie groupoids.  For instance, in \cite{pronk}, Pronk localises \'etale Lie groupoids using the bicategory of fractions outlined in Section \ref{s:localisation}, and  in \cite{roberts2012}, Roberts localises Lie groupoids, proper Lie groupoids, \'etale Lie groupoids, and \'etale proper Lie groupoids using the method of anafunctors of Section \ref{s:localisation}.

We generalise these scattered results below by considering action groupoids satisfying the properties $\calP$.  In fact, we go further.  We begin by considering action groupoids with so-called equivariant functors between them.

	\begin{definition}\label{d:eqvt}
		Let $\eqactgpd_\calP$ be the sub-$2$-category of $\LieGpoid$ whose objects are  action Lie groupoids $\calG=G\ltimes X$ satisfying $\calP$, $1$-cells are \textbf{equivariant functors} (a functor $\varphi\colon G\ltimes X\to H\ltimes Y$ is \textbf{equivariant} if there exists a Lie group homomorphism $\widetilde{\varphi}\colon G\to H$ and a smooth map $\varphi_0\colon  X \to Y$ such that the functor satisfies $$\varphi_1(g,x)=(\widetilde{\varphi}(g),\varphi_0(x))$$   for all $(g,x)\in G\times X$), and $2$-cells are natural transformations. Denote by $\we_\calP$ the class of equivariant weak equivalences and $\ssw_\calP$ the class of equivariant surjective submersive weak equivalences. 
	\end{definition}

Equivariant surjective submersive weak equivalences take on a very special form. 

	\begin{lemma}\label{l:fund ess equiv1}
        Let $G\ltimes X$ and $H\ltimes Y$ be action groupoids and let $\varphi\colon G\ltimes X\to H\ltimes Y$ be in an equivariant surjective submersive weak equivalence.  There is a closed normal subgroup $K\trianglelefteq G$ that acts freely and properly on $X$ and an equivariant isomorphism of Lie groupoids $$\psi\colon\lfaktor{G}{K}\ltimes\lfaktor{X}{K}\overset{\cong}{\longrightarrow}H\ltimes Y$$ satisfying $\varphi = \psi\circ \pi$ where $\pi\colon G\ltimes X\to\lfaktor{G}{K}\ltimes\lfaktor{X}{K}$ is the quotient functor.  Moreover, if $\G\ltimes X$ satisfies $\calP$ (excluding the property ``effective'') then $\psi$ is a morphism in $\eqactgpd_\calP$.
	\end{lemma}

	\begin{proof}
        It follows from the smooth fully faithfulness of $\varphi$ that the group homomorphism $\widetilde{\varphi}\colon G\to H$ is an epimorphism.  Let $K=\ker\widetilde{\varphi}$.  By the First Isomorphism Theorem for groups, the continuity of $\varphi$, and \cite[Corollaries 1.10.10, 1.11.5]{DK}, there is an induced Lie group isomorphism $\widetilde{\psi}\colon\lfaktor{G}{K}\to H$.

        It also follows from the smooth fully faithfulness of $\varphi$ that $K$ acts freely on $X$; these in turn imply that $K\ltimes X$ is Lie groupoid isomorphic to the submersion groupoid $X\ftimes{\varphi}{\varphi}X\toto X$ via the map $(s,t)\colon K\times X\to X\times X$.  Since submersion groupoids are proper, $K$ acts on $X$ freely and properly.  It follows from the Quotient Manifold Theorem (see, for instance, \cite[Theorem 1.11.4]{DK}) that $\lfaktor{G}{K}\ltimes\lfaktor{X}{K}$ is an action Lie groupoid whose action is induced by the residual action of $G$ on $\lfaktor{G}{K}$.

        Equivariance of $\varphi$ yields the identity $\varphi(kx)=\varphi(x)$ for all $k\in K$ and $x\in X$, from which it follows that $\psi_0([x]):=\varphi(x)$ is well-defined and smooth for all $[x]\in\lfaktor{X}{K}$.  Injectivity of $\psi_0$ follows from the smooth fully faithfulness of $\varphi$, and surjective submersivity of $\psi_0$ follows from that of $\varphi$.  The smooth functor $\psi:=(\widetilde{\psi},\psi_0)$ is a Lie groupoid isomorphism, as desired.

        By Remark~\ref{r:morita equivalent}, if $G\ltimes X$ is free, locally free, transitive, proper, or is an orbifold groupoid, then so is $H\ltimes Y$.  If $G$ is compact or discrete, then it is immediate that $H$ satisfies the same property.  It follows that if $G\ltimes X$ is in $\eqactgpd_\calP$, then so is $\psi$.
	\end{proof}

An equivariant weak equivalence in $\we_\calP$ also has a special form: it decomposes into the composition of an equivariant inclusion functor with an equivariant surjective submersive weak equivalence, which are in $\we_\calP$ and $\ssw_\calP$, resp.  This decomposition was originally observed in \cite{pronk-scull} in the proper \'etale case.  The discrete case of the decomposition already appears in \cite[Theorem 5.4]{CdHP} in the context of discrete dynamical systems.

Before proving the decomposition, we first show that the inclusion functor appearing in the decomposition is a weak equivalence.  Recall that if $K$ is a closed Lie subgroup of $G$ and $K\ltimes X$ is a Lie group action, then the anti-diagonal action of $K$ on $G\times X$ is free and proper and so the quotient $G\times_K X$ is a manifold by the Quotient Manifold Theorem \cite[Theorem 1.11.4]{DK}.

	\begin{lemma}\label{l:fund ess equiv2}
		Given a closed subgroup $K\leq G$ and an action groupoid $K\ltimes X$, the induced inclusion functor $i\colon K\ltimes X\to G\ltimes(G\times_K X)$ sending $(k,x)$ to $(k,[1_G,x])$ is an equivariant weak equivalence.  If either action groupoid is free, locally free, transitive, proper, or an orbifold groupoid, then so is the other.  Moreover, if $K\ltimes X$ is effective; then so is  $G\ltimes(G\times_K X)$; if $G$ is compact or discrete, then $K$ has the same property.
	\end{lemma}

	\begin{proof}
		The manifold $G\times_K X$ comes equipped with a $G$-action $g[g',x]:=[gg',x]$, for which $i$ is an equivariant functor with respect to the $K$- and $G$-actions.  Since the quotient map $G\times X\to G\times_K X$ is a principal $K$-bundle, it follows that $i$ is a weak equivalence.  If either groupoid is free, locally free, transitive, proper, or an orbifold groupoid, so is the other by Remark~\ref{r:morita equivalent}.

        If the $G$-action is not effective, then there is some $g'\in G$ contained in every stabiliser $\Stab_G([g,x])$ for $[g,x]\in G\times_K X$.  Since $\Stab_G([g,x])=g\Stab_K(x)g^{-1}$ for every $[g,x]\in G\times_K X$ it follows that $g'\in\Stab_K(x)$ for every $x\in X$, in which case the $K$-action is not effective. Finally, if $G$ is compact or discrete, it is immediate that $K$ also shares this property.
	\end{proof}

We now establish the decomposition of an equivariant weak equivalence.

	\begin{theorem}\label{t:decomposition}
		Let $G\ltimes X$ and $H\ltimes Y$ satisfy properties $\calP$ (except for ``effective''), and let $\varphi\colon G\ltimes X\to H\ltimes Y$ be an equivariant weak equivalence in which $\widetilde{\varphi}\colon G\to H$ is proper.  Then $\varphi$ factors as $i\circ\pi$ where $\pi\in\ssw_\calP$, and $i\in \we_\calP$ is an inclusion functor of the form as in Lemma~\ref{l:fund ess equiv2}.
	\end{theorem}
	
	\begin{proof}
		\textbf{Claim 1:} $\lfaktor{G}{\ker(\widetilde{\varphi})}$ is a Lie group isomorphic to $\im(\widetilde{\varphi})$.\\
            Indeed, $\lfaktor{G}{\ker(\widetilde{\varphi})}$ is a Lie group \cite[Corollary 1.11.5 and Proposition 1.11.8]{DK}, and $\widetilde{\varphi}$ descends to a smooth bijective homomorphism $\widehat{\varphi}\colon \lfaktor{G}{\ker(\widetilde{\varphi})}\to \im(\widetilde{\varphi})$.  Since $\widetilde{\varphi}$ is proper, it is closed, and so $\im(\widetilde{\varphi})$ is a closed subgroup of $H$, and hence a Lie subgroup \cite[Corollary 1.10.7]{DK}.  By \cite[Corollary 1.10.10]{DK}, $\widehat{\varphi}$ is a Lie group isomorphism.
		
		\noindent\textbf{Claim 2:} $\varphi_0(X)$ is an injectively immersed submanifold diffeomorphic to $\lfaktor{X}{\ker(\widetilde{\varphi})}$.\\
		It follows from the equivariance of $\varphi$ and the injectivity of $\Ff{\varphi}$ that $\ker(\widetilde{\varphi})$ acts freely on $X$.  Since $\widetilde{\varphi}$ is proper, $\ker(\widetilde{\varphi})$ is a compact submanifold of $G$.  Thus $\lfaktor{X}{\ker(\widetilde{\varphi})}$ is a manifold \cite[Theorem 1.11.4]{DK}.  Since $\varphi_0$ is $\ker(\widetilde{\varphi})$-invariant, it descends to a smooth surjection $\psi\colon\lfaktor{X}{\ker(\widetilde{\varphi})}\to\varphi_0(X)$.  Now suppose $x,x'\in X$ such that $\psi([x])=\psi([x'])$.  Then $\varphi_0(x)=\varphi_0(x')$, and since $\Ff{\varphi}$ is a diffeomorphism, there exists a (unique) $k\in\ker(\widetilde{\varphi})$ such that $x=k\cdot x'$.  It follows that $\psi$ is a smooth bijection onto $\varphi_0(X)$.
		
		Fix a curve $p=y_\tau\colon I\to \varphi_0(X)$.  Shrinking $I$, since $\ES{\varphi}$ is surjective submersive, there is a lift $q=(x_\tau,(h_\tau,y_\tau))\colon I\to X\ftimes{\varphi_0}{t}(H\times Y)$ of $p$.  By the smooth fully faithfulness of $\varphi$, the curve $h_\tau$ is contained in $\im\left(\widetilde{\varphi}\right)$.  By Claim 1, we identify $\im\left(\widetilde{\varphi}\right)$ with $\lfaktor{G}{\ker{\varphi}}$, and since $G\to\lfaktor{G}{\ker(\varphi)}$ is a principal $(\ker(\varphi))$-bundle, after shrinking $I$ again, there is a lift $g_t$ of $h_t$ to $G$.  But then $y_t=\psi\left([s(g_t)]_{\ker(\varphi)}\right)$, which proves that $\psi$ is an immersion.  This proves Claim 2.
	
		By Claim 2 and Lemma~\ref{l:fund ess equiv1}, $\pi:=(\widetilde{\varphi},\varphi_0)\colon G\ltimes X\to\lfaktor{G}{\ker(\varphi)}\ltimes\varphi_0(X)$ is an equivariant weak equivalence.  It is straightforward to check that $(\widehat{\varphi},\psi)$ is an isomorphism of Lie groupoids between $\lfaktor{G}{K}\ltimes\lfaktor{X}{K}$ and $\im\left(\widetilde{\varphi}\right)\ltimes\varphi_0(X)$; we identify these.  Let $i=(i_{\im(\widetilde{\varphi})},i_{\im(\varphi_0)})$, where the two components are inclusions of the images into $H$ and $Y$, resp.  By Claim 1, we have the following factorisation $$G\ltimes X\overset{\pi}{\longrightarrow}\im\left(\widetilde{\varphi}\right)\ltimes\varphi_0(X)\overset{i}{\longrightarrow}H\ltimes Y.$$  To obtain the desired decomposition, by Lemma~\ref{l:fund ess equiv2}, it suffices to show that $Y$ is $H$-equivariantly diffeomorphic to $H\times_{\im(\widetilde{\varphi})}\varphi_0(X)$.

		Define $\chi\colon H\times_{\im(\widetilde{\varphi})}\varphi_0(X)\to Y$ to be the smooth map given by $\chi([h,\varphi_0(x)]):=h\cdot\varphi_0(x)$.  Suppose $\chi([h,\varphi_0(x)])=\chi([h',\varphi_0(x')])$.  Since $\Ff{\varphi}$ is a diffeomorphism, there exists a unique $g\in G$ such that $\Ff{\varphi}(g,x)=(x,x',((h')^{-1}h,\varphi_0(x)))$, and so $x'=g\cdot x$ and $\widetilde{\varphi}(g)=(h')^{-1}h$.  Thus $[h,\varphi_0(x)]=[h',\varphi_0(x')]$, from which it follows that $\chi$ is injective.  For a fixed $y\in Y$, since $\ES{\varphi}$ is surjective, there exists $(x,(h,y))\in X\ftimes{\varphi}{t}(H\times Y)$ with $(h, y):  y \to \varphi(x)$ and thus $y=h^{-1}\varphi(x)=\chi(h^{-1},\varphi(x))$.  Thus $\chi$ is bijective.
		
		Let $p=y_\tau \colon I\to Y$ be a curve.  Since $\ES{\varphi}$ is a surjective submersion, after shrinking $I$, there is a lift $q=(x_\tau,(h_\tau,y_\tau))$ of $p$ to $X\ftimes{\varphi_0}{t}(H\times Y)$.  The curve $[h_\tau^{-1},\varphi(x_\tau)]$ has image $p$ via $\chi$, and thus $\chi$ is an immersion.  Since immersive bijections are diffeomorphisms, this shows that $\varphi$ decomposes into $i\circ\pi$ as desired.

		It remains to show that the domain of $i$ is an action groupoid satisfying $\calP$ (except for ``effective'').  But this follows from the preservation of these properties as stated in Lemmas~\ref{l:fund ess equiv1} and \ref{l:fund ess equiv2}.
	\end{proof}

We now use the recipe of Roberts \cite{roberts2021} to produce a localisation $\ana_\calP$ of $\eqactgpd_\calP$ at $\we_\calP$, which has the  equivariant anafunctors as $1$-cells.  Using \cite[Theorem 3.24]{roberts2021}, we show that this bicategory is equivalent to $\eqactgpd_\calP[\we_\calP^{-1}]$.  Finally, we show that both of these bicategories are equivalent to $\LieGpoid[\we^{-1}]_\calP$, the full sub-bicategory of $\LieGpoid[\we^{-1}]$ whose objects are action groupoids satisfying $\calP$.

We begin by constructing our localisation of $\eqactgpd_\calP$ at $\we_\calP$.
			
	\begin{definition}\label{d:eqvt ana morph}
	    A \textbf{$\calP$-anafunctor} is a generalised morphism $$G\ltimes X\underset{\varphi}{\ot}K\ltimes Y\underset{\psi}{\longrightarrow}H\ltimes Z$$ in which all groupoids involved are action groupoids satisfying $\calP$, with $\psi$ equivariant and $\varphi\in \ssw_\calP$. 
	\end{definition}
	
	We will compose $\calP$-anafunctors using the strict pullback.  So we need to verify that $\eqactgpd_\calP$  is also closed under strict pullbacks.    

	\begin{lemma}\label{l:eqvt strict pullback}
		Let $\calG = G \ltimes X$, $\calH = H \ltimes Y$, and $\calK = K \ltimes Z$, and let  $\varphi\colon\calG\to\calK$ and $\psi\colon\calH\to\calK$ be equivariant functors in $\eqactgpd_\calP$ (with their defining group maps $\widetilde{\varphi}, \widetilde{\psi}$) such that $\varphi\in \ssw_\calP$.  The strict pullback groupoid $\calG\ftimes{\varphi}{\psi}\calH$ is an action groupoid of a $(G\ftimes{\widetilde{\varphi}}{\widetilde{\psi}}H)$-action that satisfies properties $\calP$, and $\pr_1$ and $\pr_2$ are equivariant with respect to the restricted projection homomorphisms from $G\ftimes{\widetilde{\varphi}}{\widetilde{\psi}}H$ with $\pr_2\in \ssw_\calP$.
	\end{lemma}

	\begin{proof}
		By Proposition~\ref{p:strict pullback}, $\calG\ftimes{\varphi}{\psi}\calH$ is a Lie groupoid and $\pr_2\in \ssw$.  Since $G\ftimes{\widetilde{\varphi}}{\widetilde{\psi}}H$ is a closed subgroup of the Lie group $G\times H$, it is a Lie subgroup \cite[Corollary 1.10.7]{DK}.  It is straightforward to check that $\calG\ftimes{\varphi}{\psi}\calH$ is isomorphic to $(G\ftimes{\widetilde{\varphi}}{\widetilde{\psi}}H)\ltimes(X\ftimes{\varphi_0}{\psi_0}Y)$.   The restricted projection functors $\pr_1$ and $\pr_2$ from $\calG\ftimes{\varphi}{\psi}\calH$ are equivariant with respect to the restricted projection functors on $G\ftimes{\widetilde{\varphi}}{\widetilde{\psi}}H$.

		If $\calH$ is free, locally free, transitive, proper, or an orbifold groupoid, then so is $\calG\ftimes{\varphi}{\psi}\calH$ since $\pr_2$ is a weak equivalence, and these properties are Morita invariants; see Remark~\ref{r:morita equivalent}.  If $G$ and $H$ are compact/discrete, then so is $G\ftimes{\widetilde{\varphi}}{\widetilde{\psi}}H$.  An examination of the stabilisers of the $(G\ftimes{\widetilde{\varphi}}{\widetilde{\psi}}H)$-action on $X\ftimes{\varphi_0}{\psi_0}Y$ yields that if the $G$- and $H$-actions are effective, then so is the $(G\ftimes{\widetilde{\varphi}}{\widetilde{\psi}}H)$-action.
	\end{proof}

	\begin{remark}\label{r:eqvt strict pullback}
		Lemma~\ref{l:eqvt strict pullback} in fact does not require $\calK$ to satisfy $\calP$; this property was not used in the proof. 
	\end{remark}

Thus we know that we can define the composition of $\calP$-anafunctors using the strict pullback and get another $\calP$-anafunctor.   We can now construct a bicategory localising $\eqactgpd_\calP$ following the anafunctor method of Section~\ref{s:localisation},  provided $\ssw_\calP$ is a so-called ``bi-fully faithful singleton strict pretopology'', using the terminology of Roberts \cite[Definitions 2.9, 2.12]{roberts2021}.  We have already verified that $\ssw_\calP$ satisfies the conditions this entails:   All identity arrows are in $\ssw_\calP$, which is immediate.   $\ssw_\calP$ must be closed under strict pullback, which is  Lemma~\ref{l:eqvt strict pullback}.  $\ssw_\calP$ must be closed under composition, which follows from Item~\ref{i:3for2} of Lemma~\ref{l:weak equiv} and the fact that surjective submersions and equivariant maps are closed under composition.  Finally, elements of $\ssw_\calP$ inherit ``representably fully faithfulness'' (equivalent to smoothly fully faithfulness, which is Item~\ref{i:ff} of Lemma~\ref{l:weak equiv}) and co-fully faithfulness (Lemma~\ref{l:co-ff}) from $\ana$.  Thus, by \cite[Theorem 3.20]{roberts2021} we have:

	\begin{proposition}\label{p:eqvt bicat}
		There is a bicategory $\gana_\calP$ whose objects are those of $\eqactgpd_\calP$, arrows are $\calP$-anafunctors, and $2$-cells are the natural transformations of  Definition~\ref{d:ana gen morph}.
	\end{proposition}

To compare $\gana_\calP$ to $\eqactgpd_\calP[\we_\calP^{-1}]$ using Roberts' setup, we have to confirm that ``weak equivalences'' as defined by Roberts are the same as ours here.  Again, ``representable fully faithfulness'' is the same as smooth fully faithfulness by  Item~\ref{i:ff} of Lemma~\ref{l:weak equiv}.  We check the $\ssw_\calP$-locally split condition using the following lemma showing that $\eqactgpd_\calP$ admits weak pullbacks, from which the required $\ssw_\calP$-locally split condition follows.

	\begin{lemma}\label{l:eqvt weak pullback}
		Let $\calG=G\ltimes X$ and $\calH=H\ltimes Y$ be objects in $\eqactgpd_\calP$, and let $\varphi\colon\calG\to\calK$ and $\psi\colon\calH\to\calK$ be functors with $\varphi\in \we$.  Then the weak pullback $\calG\wtimes{\varphi}{\psi}\calH$ is isomorphic to an action groupoid of a ($G\times H$)-action on its object space $Z$ that satisfies properties $\calP$, and $\pr_1$ and $\pr_3$ are equivariant with respect to the projection homomorphisms from $G\times H$ with $\pr_3 \in \ssw_\calP$.   
	\end{lemma}

	\begin{proof}
		By Proposition~\ref{p:weak pullback}, $\calG\wtimes{\varphi}{\psi}\calH$ is a Lie groupoid and $\pr_3\in \ssw$.   Let $Z$ be its object space.  The action of $G\times H$ on $Z$ given by
$$((g,h),(x,k,y))\mapsto (gx,\psi(h,y)k\varphi(g,x)^{-1},hy)$$ yields an action groupoid canonically isomorphic as a Lie groupoid to $(G\times H)\ltimes Z$ to $\calG\wtimes{\varphi}{\psi}\calH$, and $\pr_1$ and $\pr_3$ are equivariant with respect to the projection homomorphisms $G\times H\to G$ and $G\times H\to H$ respectively.

        If $\calH$ is free, locally free, transitive, proper, or an orbifold groupoid, then so is $\calG\wtimes{\varphi}{\psi}\calH$ by Remark~\ref{r:morita equivalent}.  If $G$ and $H$ are compact/discrete, then so is $G\times H$.  Finally, an examination of the stabilisers of the ($G\times H$)-action on $Z$ immediately yields that if the $G$- and $H$- actions are effective, then so is the $(G\times H)$-action.
	\end{proof}

    That elements of $\we_\calP$ are $\ssw_\calP$-locally split now follows from  Lemma~\ref{l:loc split} and Lemma~\ref{l:eqvt weak pullback}.  It now follows from \cite[Theorem 3.24]{roberts2021} that:

	\begin{proposition}\label{p:eqvt bicat 2}
		The inclusion $\gana_\calP\to\eqactgpd_\calP[\we_\calP^{-1}]$ is an equivalence of bicategories, where this inclusion takes a $2$-cell to its equivalence class.
	\end{proposition}

\section{The Equivalence of $\gana_\calP$ and $\LieGpoid\lbrack \we^{-1}\rbrack_\calP$}\label{s:eqvt-to-full}

In the previous section, we constructed  a bicategory $\gana_\calP$ out of the action groupoids satisfying the chosen properties $\calP$, with $1$-cells given by $\calP$-anafunctors.    In this section, we prove Theorem~\ref{t:eqvt Lie gpds}, which states that $\gana_\calP$ and $\eqactgpd_\calP[\we_\calP^{-1}]$ are equivalent to the full sub-bicategory $\LieGpoid[\we^{-1}]_\calP$ of $\LieGpoid[\we^{-1}]$ whose objects are exactly the action groupoids satisfying $\calP$.  This allows us to replace any generalised morphism between two such groupoids with one from either $\gana_\calP$ or $\eqactgpd_\calP[\we_\calP^{-1}]$: these all admit $2$-cells between each other in $\LieGpoid[W^{-1}]$.

The objects of $\gana_\calP$ and $\LieGpoid[\we^{-1}]_\calP$ are the same, and every 1-cell of $\gana_\calP$ is a particular kind of generalised morphism, thus defining a 1-cell in $\LieGpoid[\we^{-1}]_\calP$.   Similarly, every 2-cell of $\gana_\calP$ also represents a 2-cell of $\LieGpoid[\we^{-1}]_\calP$.  So we have an inclusion:  

	\begin{definition}\label{d:IP}
		Define $I_\calP\colon\gana_\calP\to\LieGpoid[\we^{-1}]_\calP$ to be the assignment sending objects to themselves, sending a $\calP$-anafunctor to itself as a generalised morphism, and sending a $2$-cell between $\calP$-anafunctors to its equivalence class as a $2$-cell between generalised morphisms.
	\end{definition}
	
	The goal of this section is to show that this inclusion is a pseudofunctor which induces an equivalence of bicategories.    Thus we have to check that the inclusion $I_\calP$ respects compositions and unitors, and that it is essentially surjective and fully faithful:  any generalised morphism between two objects of $\LieGpoid[\we^{-1}]_\calP$ admits a $2$-cell from itself to a $\calP$-anafunctor, and that any $2$-cell between $\calP$-anafunctors can be represented by a unique $2$-cell from $\gana_\calP$.  

We begin with the generalised morphisms.  Our strategy will be to show that any generalised morphism is equivalent to a $\calP$-anafunctor induced by a bibundle.  The theory of bibundles offers another method of localising $\LieGpoid$ with a more geometric flavour; see \cite{haefliger,HS,lerman} for details.  We do not require the full theory here, but simply borrow the necessary concepts. 

	\begin{proposition}\label{p:gen morph equiv eqvt}
		Any generalised morphism $$\calG=G\ltimes X\underset{\varphi}{\ot}\calK\underset{\psi}{\longrightarrow}H\ltimes Z=\calH$$ between objects in $\eqactgpd_\calP$ admits a $2$-cell from itself to a $\calP$-anafunctor $\calG\underset{\chi}{\ot}\calL\underset{\omega}{\longrightarrow}\calH.$  
	\end{proposition}

	\begin{proof} Given the generalised morphism $\calG\underset{\varphi}{\ot}\calK\underset{\psi}{\longrightarrow}\calH$, one can construct an anafunctor $$\calG\underset{\chi}{\ot}\calL\underset{\omega}{\longrightarrow}\calH$$ in which $\calL$ is the action groupoid $\calG\ltimes\calL_0\rtimes\calH$ as constructed in the proof of \cite[Theorem 4.6.3]{del hoyo}.  Since $\calG=G\ltimes X$ and $\calH=H\ltimes Z$, it follows that $\calL$ is a $(G\times H)$-action groupoid. Thus all that remains is to verify that $\calL$ inherits the properties in $\calP$.  

    Since $\chi$ is a weak equivalence, $\calL$ is free, locally free, transitive, proper, or an orbifold groupoid if $\calG$ is, by Remark~\ref{r:morita equivalent}. 
 If $G$ and $H$ are compact/discrete, then so is $G\times H$.  Finally, if the actions of $G$ and $H$ on $X$ and $Y$, resp., are effective, then an examination of the stabilisers of the ($G\times H$)-action of $\calL$ reveals that this action is effective as well.
    \end{proof}

	We now want to prove a result similar to Proposition~\ref{p:gen morph equiv eqvt} for $2$-cells: that any $2$-cell between $\calP$-anafunctors in $\LieGpoid[W^{-1}]$ has a unique representative that is a $2$-cell in $\gana_\calP$. The proof of this  will require several lemmas, following the outline of the proof of a similar result of Pronk-Scull (see \cite[Section 5]{pronk-scull:bicategories}), but with some necessary modifications: equivariant surjective submersive weak equivalences are not preserved under natural transformations, so we cannot follow Pronk-Scull verbatim.  Existence of the $2$-cell representative from $\gana_\calP$ is the content of the first lemma below, which is a modified version of \cite[Lemma 5.2]{pronk-scull:bicategories}.

	\begin{lemma}\label{l:2-cell}
		Let $\calG=G\ltimes X$ and $\calH=H\ltimes Y$ be objects of $\eqactgpd_\calP$. Suppose we have two $\mathcal{P}$-anafunctors, the top and bottom of the diagram below, with a 2-cell connecting them in $\LieGpoid[\we^{-1}]_\calP$ represented by the following diagram:
   \begin{equation}\label{e:L-2cell}
		\xymatrix{
		 & & \calK \ar@{>>}[dll]_{\varphi}^{\simeq} \ar@{>>}[drr]^{\psi} & & \\
		\calG & \Downarrow \eta_1 & \calL \ar[u]_{\simeq}^{\alpha} \ar[d]^{\simeq}_{\alpha'} & \Downarrow \eta_2 & \calH \\
		 & & \calK'. \ar@{>>}[ull]^{\varphi'}_{\simeq} \ar@{>>}[urr]_{\psi'} & & \\
		}
		\end{equation}
	This $2$-cell is represented by the following $2$-cell from $\gana_\calP$: 
		\begin{equation}\label{e:2-cell}
		\xymatrix{
		 & & \calK \ar@{>>}[dll]_{\varphi}^{\simeq} \ar[drr]^{\psi} & & \\
		\calG & \circlearrowright & \calK\ftimes{\varphi}{\varphi'}\calK' \ar@{>>}[u]_{\simeq}^{\pr_1} \ar@{>>}[d]^{\simeq}_{\pr_2} & \Downarrow \nu & \calH \\
		 & & \calK'. \ar@{>>}[ull]^{\varphi'}_{\simeq} \ar[urr]_{\psi'} & & \\
		}
		\end{equation}
	\end{lemma}

	\begin{proof} 
		By Lemma~\ref{l:eqvt strict pullback}, $\calK\ftimes{\varphi}{\varphi'}\calK'$ is an action groupoid of a Lie group action of $K:=(G\times H)\ftimes{\widetilde{\varphi}}{\widetilde{\psi}}(G\times H)$.
		
	Define $\widetilde{\calL}:=(\calK\ftimes{\varphi}{\varphi'}\calK')\wtimes{\pr_1}{\alpha}\calL$, and consider the following diagram, in which Proposition~\ref{p:weak pullback} justifies the decorations on the arrows:
		$$\xymatrix{
		 & & \calK \ar@{>>}@/_1pc/[ddll]_{\varphi}^{\simeq} \ar@/^1pc/[ddrr]^{\psi} & & \\
		 & & \overset{\PR_2}{\Rightarrow} & & \\
		\calG & \calK\ftimes{\varphi}{\varphi'}\calK' \ar@{>>}@/^/[uur]^{\simeq}_{\pr_1} \ar@{>>}@/_/[ddr]_{\simeq}^{\pr_2} & \widetilde{\calL} \ar@{>>}[l]_{\quad\quad\simeq}^{\quad\quad\pr_1} \ar@{>>}[r]^{\simeq}_{\pr_3} & \calL \ar@/_/[uul]_{\simeq}^{\alpha} \ar@/^/[ddl]^{\simeq}_{\alpha'} & \calH \\
		 & & & & \\
		 & & \calK'. \ar@{>>}@/^1pc/[uull]^{\varphi'}_{\simeq} \ar@/_1pc/[uurr]_{\psi'} & & \\
		}$$
By Item~\ref{i:ff} of Lemma~\ref{l:weak equiv}, the natural transformation $$(\eta_1\pr_3)\circ(\varphi\PR_2)\circ(\ID_{\varphi'\circ\pr_2}\pr_1)\colon\varphi'\circ\pr_2\circ\pr_1\Rightarrow\varphi'\circ\alpha'\circ\pr_3$$ factors as $\varphi'\mu$ for a unique natural transformation $\mu\colon\pr_2\circ\pr_1\Rightarrow\alpha'\circ\pr_3$, making the lower triangle in the above diagram $2$-commute.  It follows from the definition of $\mu$ that $$(\eta_1\pr_3)\circ(\varphi\PR_2)=(\varphi'\mu)\circ(\ID_{\varphi\circ\pr_1}\pr_1).$$  By Lemma~\ref{l:co-ff}, the natural transformation 
$$(\psi'\mu^{-1})\circ(\eta_2\pr_3)\circ(\psi\PR_2)\colon\psi\circ\pr_1\circ\pr_1\Rightarrow\psi'\circ\pr_2\circ\pr_1$$ factors as $\nu\pr_1$ for a unique natural transformation $\nu\colon\psi\circ\pr_1\Rightarrow\psi'\circ\pr_2$.  It follows from the definition of $\nu$ that $$(\eta_2\pr_3)\circ(\psi\PR_2)=(\psi'\mu)\circ(\nu\pr_1).$$  This shows that the diagram \eqref{e:2-cell} is indeed a representative of a $2$-cell, is in $\gana_\calP$, and in the same equivalence class as the diagram \eqref{e:L-2cell}.
	\end{proof}

Thus we have shown that any $2$-cell between $\calP$-anafunctors is represented by a $2$-cell from $\gana_\calP$.  We now need to prove the uniqueness of such a representative.  
We begin with a technical lemma.    

	\begin{lemma}\label{l:subm repr}
		Given two representatives of a $2$-cell which both come from $\gana_\calP$ and have the form 
		$$\xymatrix{
		 & & \calK \ar@{>>}[dll]_{\varphi}^{\simeq} \ar[drr]^{\psi} & & \\
		\calG & \circlearrowright & \calK\ftimes{\varphi}{\varphi'}\calK' \ar@{>>}[u]_{\simeq}^{\pr_1} \ar@{>>}[d]^{\simeq}_{\pr_2} & \Downarrow \eta, \quad \Downarrow \mu & \calH \\
		 & & \calK', \ar@{>>}[ull]^{\varphi'}_{\simeq} \ar[urr]_{\psi'} & & \\
		}$$
the generalised morphism $$\calK\ftimes{\varphi}{\varphi'}\calK'\underset{\beta}{\ot}\calM\overset{\simeq}{\longrightarrow}\calK\ftimes{\varphi}{\varphi'}\calK'$$ inducing the equivalence between the two representatives (see Definition~\ref{d:equiv gen morph}) can be chosen so that $\beta_0$ is a surjective submersion.
	\end{lemma}

	\begin{proof}
		Since the two equivalences are in the same equivalence class, there exists a generalised morphism 
		$$\calK\ftimes{\varphi}{\varphi'}\calK'\underset{\alpha}{\ot}\calL\underset{\alpha'}{\overset{\simeq}{\longrightarrow}}\calK\ftimes{\varphi}{\varphi'}\calK'$$
and natural transformations $$\nu\colon\pr_1\circ\alpha\Rightarrow\pr_1\circ\alpha' \quad \text{and} \quad \nu'\colon\pr_2\circ\alpha\Rightarrow\pr_2\circ\alpha'$$
such that
		\begin{equation}\label{e:subm repr}
			(\ID_{\varphi\circ\pr_1}\alpha')\circ(\varphi\nu)=(\varphi'\nu')\circ(\ID_{\varphi\circ\pr_1}\alpha) \quad \text{and} \quad (\mu\alpha')\circ(\psi\nu)=(\psi'\nu')\circ(\eta\alpha).
		\end{equation}

Define $\calM:=(\calK\ftimes{\varphi}{\varphi'}\calK')\ftimes{\varphi\circ\pr_1}{\varphi\circ\pr_1\circ\alpha}\calL$.  By Item~\ref{i:3for2} of Lemma~\ref{l:weak equiv} and Proposition~\ref{p:strict pullback}, $\calM$ is a Lie groupoid.   It follows from Item~\ref{i:ff} of Lemma~\ref{l:weak equiv} and the fact that $\varphi\circ\pr_1=\varphi'\circ\pr_2$ that the two natural transformations between functors $\calM\to\calG$
$$(\varphi\circ \pr_1)\circ\pr_1=\varphi\circ(\pr_1\circ\alpha)\circ\pr_2\underset{\varphi\nu\pr_2}{\Longrightarrow}\varphi\circ(\pr_1\circ\alpha')\circ\pr_2 \quad \text{and}$$
$$(\varphi' \circ\pr_2) \circ \pr_1=(\varphi\circ\pr_1)\circ\pr_1 = (\varphi\circ\pr_1\circ\alpha)\circ\pr_2 = \varphi'\circ(\pr_2\circ\alpha)\circ\pr_2 \underset{\varphi'\nu'\pr_2}{\Longrightarrow} \varphi'\circ(\pr_2\circ\alpha')\circ \pr_2,$$
factor as $\varphi\omega$ and $\varphi'\omega'$ for some natural transformations $\omega\colon  \pr_1 \circ\pr_1 \Rightarrow \pr_1 \circ (\alpha' \circ \pr_2)$ and $\omega'\colon  \pr_2 \circ \pr_1 \Rightarrow \pr_2 \circ (\alpha' \circ \pr_2)$. Thus we have the following $2$-commutative diagram:
		$$
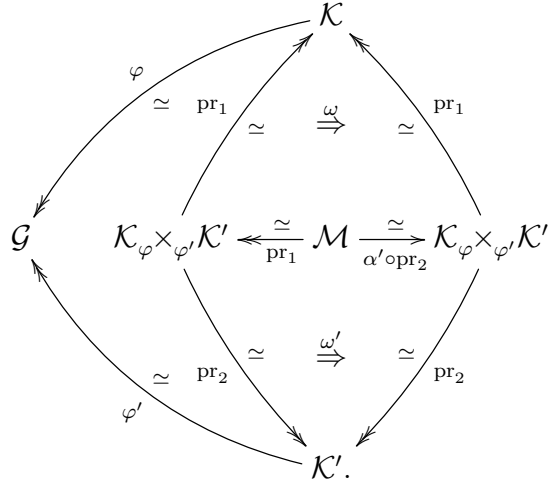
\xymatrix{
		& & \calK \ar@{>>}@/_1.5pc/[ddll]_{\varphi}^{\simeq} & \\
		& & \overset{\omega}{\Rightarrow} & \\
		\calG & \calK\ftimes{\varphi}{\varphi'}\calK' \ar@{>>}@/^/[uur]_{\simeq}^{\pr_1} \ar@{>>}@/_/[ddr]^{\simeq}_{\pr_2} & \calM \ar@{>>}[l]_{\quad\quad\simeq}^{\quad\quad\pr_1} \ar[r]^{\simeq\quad}_{\alpha'\circ\pr_2\quad} & \calK\ftimes{\varphi}{\varphi'}\calK' \ar@{>>}@/_/[uul]^{\simeq}_{\pr_1} \ar@{>>}@/^/[ddl]_{\simeq}^{\pr_2} \\
		& & \overset{\omega'}{\Rightarrow} & \\
		& & \calK'. \ar@{>>}@/^1.5pc/[uull]^{\varphi'}_{\simeq} & \\
		}$$

	To show that this diagram is an equivalence in $\LieGpoid[W^{-1}]_\calP$ between the original diagrams, we need to check Equations~\eqref{e:equiv gen morph}.  The first equation is straightforward from the definitions of $\omega$ and $\omega'$ and the first equation of \eqref{e:subm repr}.  The second follows from the smooth fully faithfulness of $\varphi$, $\varphi'$, $\pr_1\circ\alpha$, and $\pr_2\circ\alpha$; the naturality of $\nu$ and $\nu'$; Item~\ref{i:ff} of Lemma~\ref{l:weak equiv}; and the second equation of \eqref{e:subm repr}.

    Since $\varphi\circ\pr_1$ and its composition with $\alpha$ are weak equivalences, so are $\pr_1,\alpha'\circ\pr_2\colon\calM\to\calK\ftimes{\varphi}{\varphi'}\calK'$, with $(\pr_1)_0$ a surjective submersion.  Thus the generalised morphism that we require is $$\calK\ftimes{\varphi}{\varphi'}\calK'\underset{\pr_1}{\ot}\calM\underset{\alpha'\circ\pr_3}{\overset{\simeq}{\longrightarrow}}\calK\ftimes{\varphi}{\varphi'}\calK'$$ where $\beta=\pr_1$.
	\end{proof}

We now prove uniqueness of the $2$-cell representative in Lemma~\ref{l:2-cell}.

	\begin{lemma}\label{l:almost unique 2cell}
		If the two diagrams below represent the same $2$-cell in $\LieGpoid[W^{-1}]_\calP$,
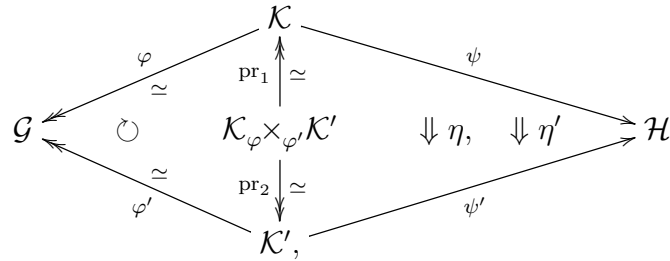
		$$\xymatrix{
		 & & \calK \ar@{>>}[dll]_{\varphi}^{\simeq} \ar[drr]^{\psi} & & \\
		\calG & \circlearrowright & \calK\ftimes{\varphi}{\varphi'}\calK' \ar@{>>}[u]_{\simeq}^{\pr_1} \ar@{>>}[d]^{\simeq}_{\pr_2} & \Downarrow \eta, \quad \Downarrow \eta' & \calH \\
		 & & \calK', \ar@{>>}[ull]^{\varphi'}_{\simeq} \ar[urr]_{\psi'} & & \\
		}
		$$
then $\eta = \eta'$.  		
	\end{lemma}

	\begin{proof}  
		Since $\eta$ and $\eta'$ are in the same equivalence class, by Lemma~\ref{l:subm repr}, there exists a generalised morphism $$\calK\ftimes{\varphi}{\varphi'}\calK'\underset{\gamma}{\ot}\calL\underset{\gamma'}{\overset{\simeq}{\longrightarrow}}\calK\ftimes{\varphi}{\varphi'}\calK'$$ and natural transformations $$\mu\colon\pr_1\circ\gamma\Rightarrow\pr_1\circ\gamma' \quad \text{and} \quad \mu'\colon\pr_2\circ\gamma\Rightarrow\pr_2\circ\gamma'$$ inducing the equivalence relation between them, in which $\gamma$ is a surjective submersive weak equivalence.  It follows from Lemma~\ref{l:co-ff} and the middle four exchange (a standard coherence relation for bicategories; see, for instance, \cite[(2.1.9)]{JY}) that $\eta\gamma=\eta'\gamma$.  Since $\gamma \in \ssw$, Lemma~\ref{l:co-ff} implies that $\eta=\eta'$, and thus the two $2$-cells are equal in $\gana$. 
	\end{proof}

	Lastly, we will verify that  the inclusion $I_\calP$ is a pseudofunctor and respects the operations in the localised bicategories.  

	\begin{proposition}\label{p:IP}
		The assignment $I_\calP\colon\gana_\calP\to\LieGpoid[\we^{-1}]_\calP$ is a pseudofunctor. 
	\end{proposition}

	\begin{proof}
		To begin, we must show that for each pair of action groupoids $\calG:= G\ltimes X$ and $\calH:= H\ltimes Y$, $I_\calP$ induces a functor $\ana_\calP(\calG,\calH)\to\LieGpoid[\we^{-1}]_\calP(\calG,\calH)$ between the categories of $1$-cells between $\calG$ and $\calH$, with $2$-cells between those.  In particular, $I_\calP$ must respect vertical composition and unit $2$-cells.

		Suppose the diagrams which are being vertically composed are in fact in $\ana_\calP$, and so $\varphi$, $\varphi'$, and $\varphi''$ are surjective submersive, $\calL_1=\calK\ftimes{\varphi}{\varphi'}\calK'$, $\calL_2=\calK'\ftimes{\varphi'}{\varphi''}\calK''$, the vertical maps from $\calL_i$ are projection maps, and $\mu_1$ and $\mu_2$ are trivial.  We consider the vertical composition in the two categories:  in $\LieGpoid[\we^{-1}]_\calP$  the vertical composition is defined by the diagram  
			$$\xymatrix{
				 & & \calK \ar@{>>}[dll]_{\varphi}^{\simeq} \ar[drr]^{\psi} & & \\
				\calG & \Downarrow\varphi'\pr_2 & (\calK\ftimes{\varphi}{\varphi'}\calK')\wtimes{\pr_2}{\pr_1}(\calK'\ftimes{\varphi'}{\varphi''}\calK'') \ar@{>>}[d]_{\pr_3\circ\pr_2}^{\simeq} \ar@{>>}[u]^{\pr_1\circ\pr_1}_{\simeq} & \Downarrow\kappa & \calH \\
				 & & \calK'' \ar@{>>}[ull]^{\varphi''}_{\simeq} \ar[urr]_{\psi''} & & \\ }$$
where $\kappa=(\eta_2\pr_3)\circ(\psi'\PR_2)\circ(\eta_1\pr_1)$ defining the vertical composition as in   \cite[Subsection 2.3]{pronk} or \cite[Section 3]{pronk-scull:bicategories}.

In $\gana_\calP$, we define the vertical compostion by 
     $$\xymatrix{
				 & & \calK \ar@{>>}[dll]_{\varphi}^{\simeq} \ar[drr]^{\psi} & & \\
				\calG & \circlearrowright & \calK\ftimes{\varphi}{\varphi''}\calK'' \ar@{>>}[d]_{\pr_2}^{\simeq} \ar@{>>}[u]^{\pr_1}_{\simeq} & \Downarrow\lambda & \calH \\
				 & & \calK'', \ar@{>>}[ull]^{\varphi''}_{\simeq} \ar[urr]_{\psi''} & & \\}$$
where  $\lambda$ is the unique natural transformation such that $\lambda\pr_{13}=(\eta_2\pr_{23})*(\eta_1\pr_{12})$ where $\pr_{ij}=(\pr_i,\pr_j)$ is the projection of $\calK\ftimes{\varphi}{\varphi'}\calK'\ftimes{\varphi'}{\varphi''}\calK''$ and $*$ denotes horizontal composition of natural transformations.  Such a $\lambda$ exists by Lemma~\ref{l:co-ff} since $\pr_{13}\in\ssw$. Thus it suffices to show for a fixed $(y,y',y'')\in\calK\ftimes{\varphi}{\varphi'}\calK'\ftimes{\varphi'}{\varphi''}\calK''$ that $\lambda\pr_{13}(y,y',y'')=\nu\pr_{13}(y,y',y'')$, where again $\nu$ is the $\nu$ of Lemma~\ref{l:2-cell}. 
 
In order to show that our inclusion respects vertical composition, we need to show that these two diagrams are in the same equivalence class.    To accomplish this, we will apply Lemma~\ref{l:2-cell} to the first diagram above, and show that the resulting natural transformation on the right (the $\nu$ of the lemma) is equal to $\lambda$.  Unravelling the vertical composition of $\LieGpoid[\we^{-1}]_\calP$ and the definition of $\nu$ in Lemma~\ref{l:2-cell}, we obtain $$\nu(y,y'')=\psi''(\mu^{-1}((\widetilde{y}'',y''),k,z))\eta_2(\widetilde{y}'_2,\widetilde{y}'')\psi'(\widetilde{k}')\eta_1(\widetilde{y},\widetilde{y}')\psi(k)$$ where $\mu$ is as defined in the proof of Lemma~\ref{l:2-cell}, $((y,y''),k,z)\in(\calK\ftimes{\varphi}{\varphi''}\calK'')\wtimes{\pr_1}{\pr_1\circ\pr_1}\calL$, and $$z:=((\widetilde{y},\widetilde{y}'_1),\widetilde{k}',(\widetilde{y}'_2,\widetilde{y}''))\in\calL:=(\calK\ftimes{\varphi}{\varphi'}\calK')\wtimes{\pr_2}{\pr_1}(\calK'\ftimes{\varphi'}{\varphi''}\calK'').$$  The result is independent of choice of (admissible) $k$ and $z$, and so we choose $z=((y,y'),u_{y'},(y',y''))$ and $k=u_y$, after which we obtain $$\nu(y,y'')=\eta_2(y',y'')\eta_1(y,y')=\lambda(y,y'').$$  Thus $I_\calP$ preserves vertical composition.

		Fix a $\calP$-anafunctor $\calG\underset{\varphi}{\ot}\calK\underset{\psi}{\longrightarrow}\calH$, and let $\Delta\colon\calK\to\calK\ftimes{\varphi}{\varphi}\calK$ be the diagonal map.  Then $\pr_1\circ\Delta=\pr_2\circ\Delta$, and $\iota_{\calG\leftarrow\calK\to\calH}\Delta$ is trivial.  Via $\calK\ftimes{\varphi}{\varphi}\calK\underset{\Delta}{\ot}\calK\overset{=}{\longrightarrow}\calK$ it follows that the identity $2$-cell of $\calG\underset{\varphi}{\ot}\calK\underset{\psi}{\longrightarrow}\calH$ in $\ana_\calP$ is equivalent to the identity $2$-cell in $\LieGpoid[\we^{-1}]_\calP$, and we conclude that $I_\calP$ induces a functor $\ana_\calP(\calG,\calH)\to\LieGpoid[\we^{-1}]_\calP(\calG,\calH)$.

		Since the identity generalised morphism of a Lie groupoid $\calG$ is the same as the identity anafunctor, $I_\calP$ trivially preserves identity $1$-cells.

		By Proposition~\ref{p:comp gen morph} for each pair of $\calP$-anafunctors   $\calG\underset{\varphi}{\ot}\calK\underset{\psi}{\longrightarrow}\calH$ and $\calH\underset{\chi}{\ot}\calL\underset{\omega}{\longrightarrow}\calI$ there is a $2$-cell in $\LieGpoid[\we^{-1}]_\calP$ from the composition as generalised morphisms to the composition as anafunctors, represented by $(\ic,\id_{\calK\ftimes{\psi}{\chi}\calL},\ID_{\varphi\circ\pr_1\circ\ic},\ID_{\pr_3\circ\omega\circ\ic})$. We now check the first of three coherence conditions (namely, \cite[(4.1.3)]{JY} or (M.1) of \cite[page 30]{benabou}), which indicates that the various compositions of $\calP$-anafunctors yields equivalent results.  Fix three  $\calP$-anafunctors
		$$\xymatrix{
			 & \calM \ar@{>>}[dl]_{\varphi}^{\simeq} \ar[dr]^{\psi} & & \calN \ar@{>>}[dl]_{\chi}^{\simeq} \ar[dr]^{\omega} & & \calP \ar@{>>}[dl]_{\xi}^{\simeq} \ar[dr]^{\zeta} & \\
			\calG & & \calH & & \calK & & \calL \\ 
		}$$
Then the first coherence condition reduces to showing that the vertical composition of the $2$-cells induced by the inclusions $\calM\ftimes{\psi}{\chi}\calN\ftimes{\omega}{\xi}\calP\to\calM\ftimes{\psi}{\chi}\calN\wtimes{\omega}{\xi}\calP$ and $\calM\ftimes{\psi}{\chi}\calN\wtimes{\omega}{\xi}\calP\to\calM\wtimes{\psi}{\chi}\calN\wtimes{\omega}{\xi}\calP$ is equal to the vertical composition of the $2$-cells induced by the inclusions $\calM\ftimes{\psi}{\chi}\calN\ftimes{\omega}{\xi}\calP\to\calM\wtimes{\psi}{\chi}\calN\ftimes{\omega}{\xi}\calP$ and $\calM\wtimes{\psi}{\chi}\calN\ftimes{\omega}{\xi}\calP\to\calM\wtimes{\psi}{\chi}\calN\wtimes{\omega}{\xi}\calP$.  These two $2$-cells are represented by the diagrams, resp.,
	$$\xymatrix{
	 & & \calM\wtimes{\psi}{\chi}\calN\wtimes{\omega}{\xi}\calP \ar@{>>}[dll]_{\varphi\circ\pr_{1}}^{\simeq} \ar[drr]^{\zeta\circ\pr_5} \\
	\calG & \circlearrowright & \calQ_1 \ar@{>>}[u]_{\simeq} \ar@{>>}[d]^{\simeq} & \circlearrowright & \calL \\
	 & & \calM\ftimes{\psi}{\chi}\calN\ftimes{\omega}{\xi}\calP \ar@{>>}[ull]^{\varphi\circ\pr_1}_{\simeq} \ar[urr]_{\zeta\circ\pr_4} \\}$$

  $$\xymatrix{
	 & & \calM\wtimes{\psi}{\chi}\calN\wtimes{\omega}{\xi}\calP \ar@{>>}[dll]_{\varphi\circ\pr_{1}}^{\simeq} \ar[drr]^{\zeta\circ\pr_{5}} \\
	\calG & \circlearrowright & \calQ_2 \ar@{>>}[u]_{\simeq} \ar@{>>}[d]^{\simeq} & \circlearrowright & \calL \\
	 & & \calM\ftimes{\psi}{\chi}\calN\ftimes{\omega}{\xi}\calP \ar@{>>}[ull]^{\varphi\circ\pr_{1}}_{\simeq} \ar[urr]_{\zeta\circ\pr_3} \\
	}$$
where $$\calQ_1:=(\calM\ftimes{\psi}{\chi}\calN\wtimes{\omega}{\xi}\calP)\wtimes{\id_{\calM\ftimes{\psi}{\chi}\calN\wtimes{\omega}{\xi}\calP}}{\id_\calM\overset{\mathrm{w}}{\times}\ic}(\calM\ftimes{\psi}{\chi}\calN\ftimes{\omega}{\xi}\calP)$$ and $$\calQ_2:=(\calM\wtimes{\psi}{\chi}\calN\ftimes{\omega}{\xi}\calP)\wtimes{\id_{\calM\wtimes{\psi}{\chi}\calN\ftimes{\omega}{\xi}\calP}}{\ic\overset{\mathrm{w}}{\times}\id_\calP}(\calM\ftimes{\psi}{\chi}\calN\ftimes{\omega}{\xi}\calP),$$ (note that we have suppressed some of the notation).  The equivalence is established by the quadruple $(j_1,j_2,\nu_1,\nu_2)$ with $\nu_1$ and $\nu_2$ trivial and where the generalised morphism $$\calQ_1\underset{j_1}{\ot}\calM\ftimes{\psi}{\chi}\calN\ftimes{\omega}{\xi}\calP\underset{j_2}{\overset{\simeq}{\longrightarrow}}\calQ_2$$ is defined by 
	\begin{align*}
		j_1(m,n,p)=&~((m,n,u_{\omega(s(n))},p),u_{\id_\calM\overset{\mathrm{w}}{\times}\ic(s(m,n,p))},(m,n,p)) \quad\text{and} \\
		j_2(m,n,p)=&~((m,u_{\chi(s(n))},n,p),u_{\ic\overset{\mathrm{w}}{\times}\id_{\calP}(s(m,n,p))},(m,n,p)).\\
	\end{align*}
Indeed, the natural transformations of Equations~\ref{e:equiv gen morph} all reduce to trivial ones.

	Finally, we check that $I_\calP$ respects the unitors from each bicategory.  For right unitors, the relevant coherence condition (see \cite[(4.1.4)]{JY} or (M.2) of \cite[page 30]{benabou})  reduces to checking that the following two diagrams represent the same $2$-cell.
		$$\xymatrix{
			& & \calK\wtimes{\psi}{\id_\calH}\calH \ar@{>>}[dll]_{\varphi\circ\pr_1}^{\simeq} \ar[drr]^{\pr_3} & & \\
			\calG & \circlearrowright & \calK\wtimes{\psi}{\id_\calH}\calH \ar@{=}[u] \ar@{>>}[d]_{\pr_1}^{\simeq} & \Uparrow\PR_2 & \calH \\
			& & \calK \ar@{>>}[ull]^{\varphi}_{\simeq} \ar[urr]_{\psi} & & \\
		}$$
		$$\xymatrix{
			& & \calK\wtimes{\psi}{\id_\calH}\calH \ar@{>>}[dll]_{\varphi\circ\pr_1}^{\simeq} \ar[drr]^{\pr_3} & & \\
			\calG & \quad\Downarrow(\varphi\circ\pr_1)\PR_2 & \calL \ar@{>>}[u]^{\ic\circ\pr_1}_{\simeq} \ar@{>>}[d]_{\pr_2\circ\pr_3}^{\simeq} & \Downarrow\rho & \calH \\
			& & \calK \ar@{>>}[ull]^{\varphi}_{\simeq} \ar[urr]_{\psi} & & \\
		}$$
where $$\calL:= (\calK\ftimes{\psi}{\id_\calH}\calH)\wtimes{\id_{\calK\times\calH}}{\pr_1}((\calK\ftimes{\psi}{\id_\calH}\calH)\ftimes{\varphi\circ\pr_1}{\varphi}\calK)$$ and $$\rho=(\rho^{\operatorname{ana}}_{\calG,\calH}(\calG\leftarrow\calK\rightarrow\calH)\pr_3)\circ(\pr_3\PR_2).$$
The equivalence is given by the quadruple $(j_1,j_2,\nu_1,\nu_2)$ where $j_1\colon\calK\to\calK\wtimes{\psi}{\id_\calH}\calH$ sends $k\in\calK_1$ to $(k,u_{s(\psi(k))},\psi(k))$ and $j_2\colon\calK\to\calL$ $k$ to $((k,\psi(k)),u_{s(k,u_{\psi(s(k))},\psi(k))},((k,\psi(k)),k))$, and both $\nu_1$ and $\nu_2$ are trivial.  Again, the natural transformations of Equations~\ref{e:equiv gen morph} are all trivial.  The computation for left unitors using the corresponding coherence condition (M.2) on page 30 of \cite{benabou} is similar.
	\end{proof}

Combining the results above yields the  desired equivalence of bicategories.

	\begin{theorem}\label{t:eqvt Lie gpds}
		The pseudofunctor $I_\calP\colon\gana_\calP\to\LieGpoid[\we^{-1}]_\calP$ is an equivalence of bicategories.  Consequently, $\gana_\calP$, $\eqactgpd_\calP[\we_\calP^{-1}]$ and $\LieGpoid[\we^{-1}]_\calP$ are all equivalent bicategories.
	\end{theorem}

	\begin{proof}
		By Proposition~\ref{p:IP}, $I_\calP$ is a pseudofunctor.  Since $I_\calP$ is surjective on objects, it suffices to show that for two action groupoids $\calG=G\ltimes X$ and $\calH=H\ltimes Y$ satisfying $\calP$, the restriction of $I_\calP$ to the category of $\calP$-anafunctors from $\calG$ to $\calH$, which maps into the category of \emph{all} generalised morphisms between them, is an equivalence of categories.  Essential surjectivity follows from Proposition~\ref{p:gen morph equiv eqvt}, and fully faithfulness follows from  Lemmas \ref{l:2-cell} and \ref{l:almost unique 2cell}.  
	\end{proof}

We  have the following immediate corollary of Theorem~\ref{t:decomposition} and Proposition~\ref{p:gen morph equiv eqvt}.

	\begin{corollary}\label{c:decomposition}
		For any subset of properties of $\calP$ which  does not include ``effective'', given a generalised morphism $\calG\underset{\varphi}{\ot}\calK\underset{\psi}{\longrightarrow}\calH$ between objects $\calG:=G\ltimes X$ and $\calH:=H\ltimes Y$ of $\eqactgpd_\calP$, there is a $2$-cell from the generalised morphism to a $\calP$-anafunctor $\calG\underset{\chi}{\ot}\calL\underset{\omega}{\longrightarrow}\calH$ in which $\chi$ is of the form $\pi$ as described in Lemma~\ref{l:fund ess equiv1}.  
	\end{corollary}

Note that if the generalised morphism in the corollary is a Morita equivalence, then $\omega$ also will decompose in the way described.  We exclude effective actions specifically from this result because they are \emph{not} preserved under equivariant weak equivalences, as the following example shows.

\begin{example}\label{x:not effective}
Let $G$ be the four-element dihedral group, with elements $(e, e), (e, \tau), (\tau, e) $ and $(\tau, \tau)$.   This acts on the set $X$ consisting of four points laid out in the cardinal directions, $N, S, E, W$:    $(\tau, e)$ reflects so that  $N$ and $S$ switch and $E, W$ are fixed,  and $(e, \tau)$ reflects so that $E$ and $W$ switch and $N, S$ are fixed, and $(\tau, \tau)$ rotates by half a turn and has no fixed points. This is an effective action.

The subgroup $K= \langle(\tau, \tau)\rangle$ acts freely, so we can take the quotient by $K$.  Then $\lfaktor{X}{K}$ consists of two points $[N]=[S]$ and $[E]=[W]$. Both of the projected points have isotropy $\mathbb{Z}/2 = \lfaktor{G}{K}$, and this action is not effective.     
\end{example}



\end{document}